\documentclass[12pt,english]{article}
\usepackage[T1]{fontenc}
\usepackage[utf8]{inputenc}
\usepackage[letterpaper]{geometry}
\usepackage{amsmath}
\usepackage{graphicx}
\usepackage{amssymb}

\makeatletter

\newtheorem{theorem}{Theorem}[section]

\newtheorem{lemma}[theorem]{Lemma}
\newtheorem{proposition}[theorem]{Proposition}

\newtheorem{remark}[theorem]{Remark}
\numberwithin{equation}{section}
\def\th{\theta}

\makeatother

\newcommand{\eqnsection}{
\renewcommand{\theequation}{\thesection.\arabic{equation}}
 \makeatletter   \csname  @addtoreset\endcsname{equation}{section}
   \makeatother}

\def\proof{\noindent\textbf{Proof: $\,$} }
\def\qed{\hfill$\Box $}

\def\({\left(}\def\){\right)}
\def\R{\mathbb{R}}

\def\E{\mathbb{E}}
\def\P{\mathbb{P}}
\def\0{\mathbf{0}}
\def\1{\mathbf{1}}

\def\de{\delta}

\def\ol{\overline}

\def\Var{{\mathop {{\rm Var\, }}}}

\makeatother

\usepackage{babel}

\makeatother

\begin{document}

\title{Annealed asymptotics for Brownian motion of renormalized
potential in mobile random medium}
\author{Xia Chen
\thanks{Supported in part by the Simons Foundation \#244767.
} \, and  \, Jie Xiong
\thanks{Supported in part by NSF grant DMS-0906907 and FDCT 076/2012/A3.}}

\maketitle

\begin{abstract}
Motivated by the study of the directed polymer model with mobile
Poissonian traps or catalysts and the stochastic parabolic Anderson model with
time dependent potential, we investigate the asymptotic behavior of
\[\E\otimes\E_0\exp\left\{\pm\theta\int^t_0\bar{V}(s,B_s)ds\right\}\qquad
(t\to\infty)\] where $\th>0$ is a constant, $\ol{V}$ is the renormalized Poisson potential
of the form
\[\ol{V}(s,x)=\int_{\R^d}\frac{1}{|y-x|^p}\(\omega_s(dy)-dy\),\]
and $\omega_s$ is the measure-valued process consisting of
independent Brownian particles whose initial positions form a
Poisson random measure on $\R^d$ with Lebesgue measure as its
intensity. Different scaling limits are obtained according to the
parameter $p$ and dimension $d$. For the logarithm of the negative exponential moment,
the range of $\frac{d}{2}<p<d$ is divided into 5 regions with
various scaling rates of the orders $t^{d/p}$, $t^{3/2}$,
$t^{(4-d-2p)/2}$, $t\log t$ and $t$, respectively. For the positive
exponential moment, the limiting behavior is studied according to
the parameters $p$ and $d$ in three regions. In the sub-critical region
($p<2$), the double logarithm of the exponential moment has a rate
of $t$. In the critical region ($p=2$), it has different behavior
over two parts decided according to the comparison of $\theta$ with
the best constant in the Hardy inequality. In the super-critical region $(p>2)$,
 the exponential moments become infinite for all $t>0$.

\begin{quote} {\footnotesize

}

\end{quote}

\end{abstract}


\renewcommand{\theequation}{\thesection.\arabic{equation}}

\section{Introduction } \label{intro}

The model of Brownian motion in a static random medium has been
thoroughly investigated. For the directed polymer with  immobile
Poissonian traps or catalysts, we refer to  Sznitman's book \cite{Sznitman-2}
for general
collection (up to the year 1998), \cite{Komorowski} for a survey and
\cite{BBH}, \cite{BTV}, \cite{GHK}, \cite{Povel} for specific
topics. For the results motivated by the parabolic Anderson models
with time-independent random potentials, we cite  \cite{CM},
\cite{CV},
 \cite{DM}, \cite{DV}, \cite{FV}, \cite{GK}, \cite{GKM},
\cite{Stolz} as a partial list of the publications on this subject.
In addition, we also point out the
very recent papers \cite{Chen-1}, \cite{CK}, \cite{CK-1} \cite{CR}
and \cite{Fu}
for the investigation related to the topic of this paper.

The case of mobile
environment is much less understood.
In the models considered in this paper,
the environment consists of  moving particles
initially distributed in the space $\R^d$ according to
a Poisson field $\omega_0(dx)$ with the Lebesgue
measure $dx$ as its intensity. The particles  move
in $\R^d$ independently so the spatial distribution of the obstacles
at the time $t$ is of the form
\begin{align}\label{poisson}
\omega_t(dx)=\int_{\R^d}\delta_{X_y(t)}(dx)\omega_0(dy)
\end{align}
where the stochastic flow $\{X_y(t)\}$ represents the paths of the
random obstacles with $X_y(0)=y$. Unless stated otherwise,
throughout this paper, $X_y(t)$ are independent Brownian motions
with the variance $\sigma^2=Var(X_y(1))>0$. More precisely, the processes
$X_y(t)-y$ are i.i.d., independent of $\omega_0(dx)$ and distributed
the same as the $\sigma$-multiple of a standard $d$-dimensional
Brownian motion.
 The notations ``$\P$'' and ``$\E$'' are used
for the
probability law and the expectation, respectively, generated by the
measure-valued process
$\{\omega_t(dx); \hskip.1in t\ge 0\}$.

Let $K(x)\ge 0$
be a properly chosen function (known as shape function) on $\R^d$.
The random field
\begin{align}\label{poisson-0}
V(t, x)=\int_{\R^d}K(y-x)\omega_t(dy)
\end{align}
represents the total potential
at $x\in\R^d$ generated by the Poisson obstacles.

In the model of directed polymer of the potential $V(t,x)$, all
possible trajectories of the random path $\{B_s\}_{0\le s\le t}$ are
re-weighted by some properly defined Gibbs measure (see
(\ref{poisson-1}) and (\ref{poisson-2}) below) that favors
trajectories in the time-space region less populated with random
traps $\{X_y(t)\}$. In this paper, $\{B_t\}_{t\ge 0}$ is a
$d$-dimensional Brownian motion independent of the environment
$\{\omega_t(dx);\hskip.05in t\ge 0\}$. Throughout this article, the
notations ``$\P_x$'' and ``$\E_x$'' are for the probability law and
the expectation, respectively, of the Brownian motion $B_s$ with
$B_0=x$.

In the quenched setting, where the set-up
is conditioned on the random environment created by the
stochastic flow $\{X_y(t)\}$,  the survival path is selected by
 the Gibbs measure
\begin{align}\label{poisson-1}
{d\mu_{t,\omega}\over d\P_0}={1\over Z_{t,\omega}}\exp\bigg\{-\theta
\int_0^tV(s, B_s)ds\bigg\}
\end{align}
defined on the space $C\([0, t];\hskip.05in\R^d\)$ of the continuous
functions $f$: $[0,t]\longrightarrow\R^d$.

In the annealed setting, where the model averages on both the Brownian
motion and the environment, the Gibbs measure is given as
\begin{align}\label{poisson-2}
{d\mu_{t}\over d(\P_0\otimes \P)}={1\over Z_{t}}\exp\bigg\{-\theta
\int_0^tV(s, B_s)ds\bigg\}.
\end{align}

In (\ref{poisson-1}) and (\ref{poisson-2}), the integral
$$
\int_0^tV(s, B_s)ds
$$
represents the accumulated potential of the Brownian path
$\{B_s\}_{0\le s\le t}$ with respect to the moving traps described
by $\omega_t(dx)$. Under the law $\mu_{t,\omega}$ or $\mu_t$,
therefore, the Brownian trajectories heavily impacted by the
Poisson obstacles are penalized and become less likely.

The Brownian motion in mobile random medium corresponds
to the parabolic Anderson model
\begin{align}\label{poisson-4}
\left\{\begin{array}{ll}\partial_tu(t,x)
=\kappa\Delta u(t,x)+\xi(t, x)u(t,x)\\\\
u(0, x)=1,\hskip.2in x\in\R^d,\end{array}\right.
\end{align}
with the time-dependent potential $\xi(t,x)$.
Take $\xi(t,x)=-\theta V(t,x)$
 and consider
two types of particles $A$ and $B$ in the space $\R^d$. The
$B$-particles evolve according to the measure-valued process
$\omega_t(dx)$. The $A$-particles diffuse from region of high
concentration to the region of low concentration according to Fick's
law, and is destroyed by nearby $B$-particles at the annihilation
rate $\theta V(t,x)$ which measures the total potential generated by
$B$-particles at $(t,x)$. The solution $u(t,x)$ of (\ref{poisson-4})
represents the time-space density of $A$-particles which are
uniformly distributed at the beginning. The story is called the
model of trapping reactions and the annihilation mechanism described
here is labeled as ``$A+B\longrightarrow B$'' in the physical
literature (see e.g., \cite{MOBC-2}),

The case when $\xi(t,x)=\theta V(t,x)$ corresponds to the branching
Brownian motion in the catalytic medium.  The $A$-particles (reactant)
diffuse according to Fick's law while under the stimulation of the
nearby $B$-particles (catalyst), each of them split into two at the rate
$\theta V(t,x)$. In this model, $u(t,x)$ represents the time-space
density of $A$-particles. The interested reader is referred to the
survey \cite{GHM} for detail physical background of this model.

The existing literature (see, e.g., \cite{DGRS},  \cite{GHM},
\cite{MOBC-2}) on the models (\ref{poisson-1}), (\ref{poisson-2})
and (\ref{poisson-4}) mainly consider the setting of the lattice
space (instead of $\R^d$) combined with pure jump random walks
(instead of Brownian motions). In these publications, the shape
function is $K(x)=\delta_0(x)$ (Dirac function). To the model of
directed polymer, this means that the only traps or catalysts that located on the
path $\{B_s\}_{0\le s\le t}$ contribute to the re-shape of the path
$\{B_s\}_{0\le s\le t}$. To the parabolic Anderson model with
$\xi(t,x)=\pm\theta V(t,x)$, it means that the $A$-particles and
$B$-particles react only in the case of collision.

Motivated by Newton's law of the universal attraction, a
renormalized Poisson potential, formally written as
\begin{align}\label{poisson-5}
\ol{V}(x)=\int_{\R^d}{1\over\vert
y-x\vert^p}\big[\omega_0(dy)-dy\big], \hskip.2in x\in\R^d,
\end{align}
is introduced in \cite{CK} in the case of the static medium (or, the case when
 $\sigma =0$).
Among other things, it has been shown (Theorem 1.1, \cite{CK})
that $\ol{V}(x)$ is well defined
as a random field if and only if $d/2<p<d$ and in this case the annealed
moment
\begin{align}\label{poisson-6}
\E\otimes\E_0\exp\bigg\{-\theta\int_0^t\ol{V}(B_s)ds\bigg\}<\infty
\end{align}
for every $\theta>0$ and $t>0$.

We intend to apply this idea to the case of mobile random medium by
constructing the renormalized potential
\begin{align}\label{poisson-7}
\ol{V}(t,x)=\int_{\R^d}{1\over\vert y-x\vert^p}\big[\omega_t(dy)-dy\big]
\hskip.2in t\ge 0,\hskip.1in x\in\R^d
\end{align}
as the replacement of $V(t,x)$ in (\ref{poisson-1}),
(\ref{poisson-2}) and in the models of the trapping reactions and
the branching  Brownian motion in catalytic media. Indeed, by a calculation of characteristic function (see Proposition 1.3 in van den Berg et al \cite{vmw}) we have that for any $t>0$
\begin{align}\label{poisson-8}
\omega_t(dx)\buildrel d\over =\omega_0(dx),
\end{align}
and hence, the field $\ol{V}(t,\cdot)$ is well defined and has the same
distribution as $\ol{V}(0,\cdot)$. Further, a slight modification of
the argument for Proposition 2.8, \cite{CK} shows that $\ol{V}(t,x)$
is continuous in probability as the function of $(t,x)$.
Consequently (Chapter 6, \cite{Borkar}), the family of random
variables (defined as equivalent classes)
$\Big\{\ol{V}(t,x);\hskip.1in (t,x)\in\R^+\times\R^d\Big\}$ yields a
measurable (joint in $(t,x)$) modification. In the remaining of the
paper, the same notation $\ol{V}(t,x)$ is used for the measurable
modification. Taking $K(x)=\vert x\vert^{-p}$ in Lemma
\ref{moment-4} below, the integral
$$
\int_0^t\ol{V}(s, B_s)ds
$$
converges almost surely under $d/2<p<d$.

The necessity of renormalization comes from the fact that
\begin{equation}
\int_{\R^d}{1\over\vert y-x\vert^p}\omega_t(dy)=\infty\hskip.2in
a.s.\hskip.1in (t,x)\in\R^+\times\R^d
\end{equation}
as $p\le d$. Without the renormalization, the Gibbs measure will be indeterminant which is $\frac00$ for trapping model and $\frac\infty\infty$ for branching model.

The rationale of renormalization is based on the
symbolic decomposition
$$
\ol{V}(t,x)=\int_{\R^d}{1\over\vert y-x\vert^p}\omega_t(dy)
-\int_{\R^d}{1\over\vert y-x\vert^p}dy
=\int_{\R^d}{1\over\vert y-x\vert^p}\omega_t(dy)
-\int_{\R^d}{1\over\vert y\vert^p}dy.
$$
The Lebesgue integral part (if treated as an ordinary constant) does
not change the mechanism of the Gibbs measures in (\ref{poisson-1})
and (\ref{poisson-2}) as it is ``renormalized'' into the partition
functions $Z_{t,\omega}$ and $Z_t$. Thus, the directed polymer of
the renormalized potential $\ol{V}(t,x)$ can be modeled in terms of
the quenched Gibbs measure
 \begin{align}\label{poisson-8'}
{d\ol{\mu}_{t,\omega}\over d\P_0}={1\over \ol{Z}_{t,\omega}}
\exp\bigg\{-\theta\int_0^t\ol{V}(s, B_s)ds\bigg\}
\end{align}
and the annealed Gibbs measure
\begin{align}\label{poisson-9}
{d\ol{\mu}_{t}\over d(\P_0\otimes\P)}={1\over \ol{Z}_{t}}
\exp\bigg\{-\theta\int_0^t\ol{V}(s, B_s)ds\bigg\}
\end{align}
in the case when the partition functions $\ol{Z}_{t,\omega}$ and
$\ol{Z}_{t}$ are finite, the condition taken neither for granted
 (given the non-locality
and singularity of the shape function $K(x)=\vert x\vert^{-p}$), nor
as the corollary to (\ref{poisson-6}).

Taking $\xi(t,x)=-\theta\ol{V}(t,x)$ in the parabolic Anderson model
(\ref{poisson-4}), we obtained a renormalized version of the model
of trapping reactions where the $B$-particles are annihilated at the
rate proportional to a Newton-type-potential and are generated at a
constant rate. With $\xi(t,x)=\theta\ol{V}(t,x)$, the branching
Brownian motions (reactants) split at the time-space rate the same
as the one in the model of trapping reactions.

For a sufficiently smooth and well-bounded function $\xi(t,x)$, by
Feynman-Kac formula (cf. Theorem 3.2 in Durrett \cite{durrett} p138)
the equation (\ref{poisson-4}) is solved by the function
\begin{align}\label{FK}
u(t,x)=\E_x\exp\bigg\{\int_0^t\xi(t-s, B_{2\kappa s})ds\bigg\}.
\end{align}
For the case when $\xi(t,x)=\pm\theta\ol{V}(t,x)$, it follows from the same
argument as in \cite{CK-1} that $u(t,x)$ is  a mild solution to
(\ref{poisson-4}) (see (\ref{th-15}) below).
 Therefore, the investigation of the random
field $u(t,x)$ defined in (\ref{FK}) is closely associated to our
understanding of the parabolic Anderson model (\ref{poisson-4}).

This paper is to study the partition function $\ol{Z}_t$ in the
annealed directed polymer defined in (\ref{poisson-9}) and the
expectation of the Feynman-Kac exponential moment given in
(\ref{FK}) with $\xi(t,x)=\pm\theta\ol{V}(t,x)$. By the identities
in law:
\begin{align}\label{poisson-10}
\ol{V}(t,x)\buildrel d\over =\ol{V}(t,0)
\end{align}
and
\begin{align}\label{poisson-11}
\Big\{\omega_{t-s}(dx);\hskip.1in 0\le s\le t\Big\}\buildrel d\over
= \Big\{\omega_{s}(dx);\hskip.1in 0\le s\le t\Big\},
\end{align}
we derive that for any $x\in\R^d$, $\theta>0$ and $t>0$,
\begin{align}\label{poisson-12}
&\E_x\exp\bigg\{\pm\theta\int_0^t\ol{V}(t-s, B_{2\kappa s})ds\bigg\}
\buildrel d\over
=\E_0\exp\bigg\{\pm\theta\int_0^t\ol{V}(t-s, B_{2\kappa s})ds\bigg\}\\
&\buildrel d\over =\E_0\exp\bigg\{\pm\theta\int_0^t\ol{V}(s,
B_{2\kappa s})ds\bigg\}=\E_0\exp\bigg\{\pm{\theta\over 2\kappa}
\int_0^{2\kappa t}\ol{V}((2\kappa)^{-1}s, B_{s})ds\bigg\}.\nonumber
\end{align}
To simplify our notation we first consider the case when $\kappa =1/2$.
Therefore, our objective is summarized into the investigation of the
annealed exponential moments
\begin{align}\label{poisson-13}
\E\otimes\E_0\exp\bigg\{\pm\theta\int_0^t\ol{V}(s,
B_{s})ds\bigg\}
\end{align}
with the concerns on integrability and long term asymptotics.

The existing literature in the setting of mobile medium deals with
the model on the lattice with the Dirac function as the shape
function (i.e., $K(x)=\delta_0(x)$). In the setting of trapping
reactions, the annealed moment
$$
\E\otimes\E_0\exp\bigg\{-\theta\int_0^tV(s, B_{s})ds\bigg\}
$$
is called the surviving probability in Bramson and Lebowitz
(\cite{BL-1} and \cite{BL-2}), Drewitz $et\;al$ \cite{DGRS}, where
it is proved that
\begin{align}\label{poisson-14}
\log\E\otimes\E_0\exp\bigg\{-\theta\int_0^tV(s, B_{s})ds\bigg\}
\sim\ \begin{cases}-c_1(\theta)\sqrt{t} & d=1 \medskip \\
-c_2(\theta){t\over\log t} & d=2 \medskip \\
-c_3(\theta)t & d\ge 3 \end{cases}
\end{align}
where $0<c_i(\theta)<\infty$ are constants and the $c_2(\theta)$ and
$c_3(\theta)$ are explicitly identified. Similar results for a
related model is also obtained by G\"artner $et\;al$ \cite{GdenHM}. Here the notation $a_t\sim b_t$ means that
$a_t/b_t\to 1$ as $t\to\infty$.

As for the branching random walks in catalytic medium, G\"artner and
den Hollander (\cite{GH}) observe a double exponential growth given
as
\begin{align}\label{poisson-15}
\lim_{t\to\infty}{1\over t}\log\log\E\otimes
\E_0\exp\bigg\{\theta\int_0^tV(s, B_{s})ds\bigg\}=C(\theta)
\end{align}
with an extended constant $0\le C(\theta)\le\infty$ (see, Theorem
1.4, \cite{GH} for the discussion on the limit $C(\theta)$). We also
refer the paper \cite{GHM-1} for the  branching random walks with
the voter model as the catalyst.

Another relevance in literature is the recent study (\cite{Chen-1},
\cite{CK}, \cite{CK-1} and \cite{CR}) on the Brownian motion of the
renormalized potential $\ol{V}(x)$ (defined in (\ref{poisson-5})) in
a static Poisson medium. It is shown (\cite{CK-1}) that under
$d/2<p<d$,
\begin{align}\label{poisson-16}
\lim_{t\to\infty}t^{-d/p}
\log\E_0\otimes\E\exp\bigg\{-\theta\int_0^t\ol{V}(B_s)ds\bigg\}=
\theta^{d/p}\omega_d{p\over d-p}\Gamma\Big({2p-d\over p}\Big)
\end{align}
for all $\theta>0$, and (\cite{CK}) that
\begin{align}\label{poisson-17}
\E_0\otimes\E\exp\bigg\{\theta\int_0^t\ol{V}(B_s)ds\bigg\}=\infty
\end{align}
for all $\theta>0$ and $t>0$, where $\omega_d$ is the volume of the $d$-dimensional unit
ball. See \cite{Chen-1} and \cite{CR} for the
investigation on the quenched setting. In view of the stationarity
(\ref{poisson-8})
of the measure-valued process $\omega_t(dx)$, we have that
$\ol{V}(t,\cdot)\buildrel d\over =\ol{V}(\cdot)$ for each $t\ge 0$.
A natural question is to ask
what will remain and what will change when the static Poisson $\omega_0(dx)$
is replaced by the mobile medium $\{\omega_t(dx);\hskip.1in t\ge 0\}$.

Finally, we point out that $V(t,x)$ is well-defined for the case of $p>d$ (we refer the reader to Proposition 2.1 of \cite{CK} for a proof of this fact). In this case, renormalization is not needed.
Because of the length of the paper, we will treat the case of $p>d$ elsewhere.

\section{Main theorems}\label{th}

We adopt all notations introduced in the previous section: the
measure-valued process $\omega_t(dx)$ is defined in (\ref{poisson})
where $\{X_y(t)\}$ are the independent Brownian motions with the
covariance matrix $\sigma I_d$ and initial location
$X_y(0)=y$, where $\sigma^2>0$ and $I_d$ is the $d\times d$ identical
matrix. The renormalized
potential $\ol{V}(t,x)$ is well defined in (\ref{poisson-7}) under
the condition
\begin{align}\label{th-1}
d/2<p<d
\end{align}

Two functions $\psi(a)$ and $\Psi(a)$ on $\R^+$
frequently appearing in the remaining
of this paper are
\begin{align}\label{moment-1}
\psi(a)=e^{-a}-1+a\hskip.1in\hbox{and}\hskip.1in
\Psi(a)=e^a-1-a,\hskip.2in a\ge 0.
\end{align}
It is easy to see that $\psi(a)$ and $\Psi(a)$ are non-negative, increasing
and convex on $\R^+$ with $\psi(a)\le \Psi(a)$.

For the negative exponential moment asymptotic behavior, the
condition (\ref{th-1}) is divided into five disjoint regimes
attached with the constants $\rho_i(\theta, \sigma^2)$ ($1\le i\le 5$),
respectively.

{\bf Regime I.} $p<2$. Write
$$
\rho_1(\theta, \sigma^2)=\theta^{d/p}\omega_d{p\over
d-p}\Gamma\Big({2p-d\over p}\Big).
$$

{\bf Regime II.} $p=2$. By (\ref{th-1}) $d=3$ in this case.
Write
$$
\rho_2(\theta,\sigma^2)=\int_{\R^3}
\E\psi\bigg(\theta\int_0^1{ds\over\vert x+X_0(s)\vert^2}\bigg)dx.
$$
The exact value of $\rho_2(\theta,\sigma^2)$
remains unknown but can be explicitly bounded (see Lemma \ref{lem1001a} below).

{\bf Regime III.}  $2<p<{d+2\over 2}$. Write
$$
\rho_3(\theta, \sigma^2)= {2^{2+d-2p\over 2}\theta^2d\omega_d\over
(2+d-2p)(4+d-2p)\sigma^{2p-d}} {\displaystyle \Gamma^2\Big({d-p\over
2}\Big) \Gamma\Big({2p-d\over 2}\Big)\over\displaystyle
\Gamma^2\Big({p\over 2}\Big)}.
$$

{\bf Regime IV.}  $p={d+2\over 2}>2$. Write
$$
\rho_4(\theta, \sigma^2)=2^{d+4\over2}d\omega_d
\Big({\theta\over (d-2)\sigma}\Big)^2.
$$

{\bf Regime V.}  $p>\max\Big\{2, {d+2\over 2}\Big\}$. The constant
$0<\rho_5(\theta, \sigma^2)<\infty$ remains unknown but
is explicitly bounded (see
(\ref{eq0413a}) below).

The next theorem gives the negative exponential moment asymptotic
behavior.

\begin{theorem}\label{th-3} Under (\ref{th-1}), we have
\begin{align}\label{th-4}
\E_0\otimes\E\exp\bigg\{-\theta\int_0^t\ol{V}(s,
B_s)ds\bigg\}<\infty \hskip.2in \theta>0,\hskip.1in t>0.
\end{align}
Further,
\begin{align}\label{th-5}
\log\E_0\otimes\E\exp\bigg\{-\theta\int_0^t\ol{V}(s, B_s)ds\bigg\}
\sim\ \begin{cases}\rho_1(\theta, \sigma^2)t^{d/p}
& \hbox{in Regime I}\medskip \\
\rho_2(\theta, \sigma^2)t^{3/2}& \hbox{in Regime II}\medskip\\
\rho_3(\theta, \sigma^2)t^{4+d-2p\over 2}& \hbox{in Regime III}\medskip \\
\rho_4(\theta, \sigma^2)t\log t& \hbox{in Regime IV}\medskip \\
\rho_5(\theta, \sigma^2)t& \hbox{in Regime V}\end{cases}
\end{align}
as $t\to\infty$.
\end{theorem}

Recall that the random potential $\ol{V}(x)$ is defined in (\ref{poisson-5}).
The comparison between (\ref{poisson-16}) and (\ref{th-5}) shows that moving
from static environment to mobile environment, the asymptotic behavior remains
the same only in the regime I (which contains the cases of $d=1$ or 2). Further, notice that counting from regimes
I to V, the deviation scales
$$
t^{d/p},\hskip.05in t^{3/2},\hskip.05in t^{4+d-2p\over 2},\hskip.05in t\log t\hskip.05in
\hbox{and}
 \hskip.05in t
$$
decrease in the following sense: For $p_1$, $p_3$ in region I, III, respectively, and for $t$ large,
\[t^{d/p_1}>t^{3/2}>t^{4+d-2p_3\over 2}>t\log t>t.\] Since the deviation scale is positively related to the survival
probability,  it is generally harder for the
Brownian particle $B_s$ to avoid the Poisson traps in  the mobile medium
than in the static medium. The fact that the deviation
scale decreases in $p$ suggests that
 the total impact of the obstacles over $B_s$
increases when the tail of the shape function gets
heavier. This observation provides the evidence showing that
the major contribution comes from the vast number of the traps located
a distance away from the Brownian particle, rather than a few
in a close neighborhood of the Brownian particle.

We now move to the case of catalytic medium.
Let $W^{1,2}(\R^d)$ be the Sobolev space defined as
$$
W^{1,2}(\R^d)=\Big\{f\in{\cal L}^2(\R^d); \hskip.1in\nabla f\in{\cal
L}^2(\R^d)\Big\}.
$$
By Lemma 7.2, \cite{Chen-1}, for any $d/2<p<\min\{2,d\}$
there is a constant $C>0$ such that
\begin{align}\label{th-6}
\int_{\R^d}{f^2(x)\over\vert x\vert^2}dx\le C\|f\|_2^{2-p}\|\nabla
f\|_2^p \hskip.2in f\in W^{1,2}(\R^d).
\end{align}
Let $\gamma(d,p)>0$ be the best constant in the above inequality.

\begin{theorem}\label{th-7} Assume (\ref{th-1}).

(1) When $p<2$,  for any $\theta>0$ and $t>0$, we have
\begin{align}\label{th-8}
\E_0\otimes\E\exp\bigg\{\theta\int_0^t\ol{V}(s, B_s)ds\bigg\}
<\infty.
\end{align}
Further,
\begin{align}\label{th-9}
\lim_{t\to\infty}{1\over t}\log\log
\E_0\otimes\E\exp\bigg\{\theta\int_0^t\ol{V}(s, B_s)ds\bigg\}
={2-p\over 2}\Big({p\over\sigma^2}\Big)^{p\over 2-p}
\Big(\theta\gamma(d,p)\Big)^{2\over 2-p}.
\end{align}

(2). When $p=2$ (and therefore $d=3$ by (\ref{th-1})), for any $t>0$,
\begin{align}\label{th-10}
\E_0\otimes\E\exp\bigg\{\theta\int_0^t\ol{V}(s, B_s)ds\bigg\}
\ \begin{cases}<\infty & \hbox{when $\theta <{\sigma^2\over 8}$}\medskip \\
=\infty & \hbox{when $\theta >{\sigma^2\over 8}$}.\end{cases}
\end{align}

Further, for any $\theta<\sigma^2/8$, we have
\begin{align}\label{th-11}
\lim_{t\to\infty}t^{-3/2}\log
\E_0\otimes\E\exp\bigg\{\theta\int_0^t\ol{V}(s, B_s)ds\bigg\}
=\int_{\R^3}\E\Psi\bigg(\theta\int_0^1{ds\over \vert
x+X_0(s)\vert^2}\bigg)dx,
\end{align}
and the right hand side is finite.

(3). When $p>2$, for any $\theta>0$ and $t>0$, we have
\begin{align}\label{th-12}
\E_0\otimes\E\exp\bigg\{\theta\int_0^t\ol{V}(s, B_s)ds\bigg\}
=\infty.
\end{align}
\end{theorem}

The comparison between (\ref{poisson-17}) and Theorem \ref{th-7}
shows that in the catalytic medium,
the mobile Poisson particles are less attractable to
the Brownian particles than static Poisson particles.
Another interesting observation is on the role
played by the parameter $p$. It is evident that the singularity of the
shape function $K(x)=\vert x\vert^{-p}$ at $x=0$ is the only reason for
the positive exponential moment to blow-up. Consequently, the model
is more likely to blow-up for greater $p$. As for the large-$t$
behavior, the singularity of the shape function is no longer the sole
driving force, as indicated by the comparison of the double exponential
growth in the sub-critical case $p<2$ versus the single exponential
growth in the critical case $p=2$.

In connection to the parabolic Anderson model posted in (\ref{poisson-4}),
write
\begin{align}\label{th-13}
u_{\pm}(t,x)=\E_x\exp\bigg\{\pm\theta\int_0^t\ol{V}(t-s, B_{2\kappa
s})ds\bigg\}.
\end{align}
By Theorems \ref{th-3}-\ref{th-7} (with $\sigma^2$ be replaced by
$(2\kappa)^{-1}\sigma^2$) and the relation (\ref{poisson-12}), $\E
u_{-}(t,x)<\infty$ and $\E u_{+}(t,x)<\infty$ for all $\theta>0$ and
$t>0$ under $d/2<p<d$ and $d/2<p<\min\{2, d\}$, respectively. In the
case $p=2$ and $d=3$, $\E u_{+}(t,x)<\infty$ if
$\theta<{\sigma^2\over 16\kappa}$; and $=\infty$ if
$\theta>{\sigma^2\over 16\kappa}$. Finally, $\E u_{+}(t,x)=\infty$
for all $\theta>0$ and $t>0$ when $p>2$.

Unfortunately, given the fact (derived from Proposition 2.9 in
\cite{CK} and the relation (\ref{poisson-8})) that for each $t>0$,
the random field $\ol{V}(t,\cdot)$ is unbounded and therefore
discontinuous in any neighborhood in $\R^d$ with positive
probability, it is unlikely that with $\xi(t,x)=\pm\theta
\ol{V}(t,x)$, the equation (\ref{poisson-4}) has path-wise solution.
On the other hand, by an argument the same as in the proofs of
Proposition 1.2 and Proposition 1.6, \cite{CK}, we can show that
whenever $u_\pm(t,x)<\infty$ a.s., particularly when $\E
u_\pm(t,x)<\infty$, $u_\pm(t,x)$ is
 the mild solution to (\ref{poisson-4})
(with $\xi(t,x)=\pm\theta\ol{V}(t,x)$), in the sense that
\begin{align}\label{th-14}
\int_0^t\int_{\R^d}p_{2\kappa(t-s)}(x-y)|\ol{V}(s,y)u_{\pm}(s,y)|\, dyds<+\infty,
\quad x\in\R^d, \, t>0
\end{align}
and
\begin{align}\label{th-15}
u_{\pm}(t,x)=1 \pm\theta\int_0^t\int_{\R^d}p_{2\kappa(t-s)}(x-y)
\ol{V}(s, y)u_\pm(s,y)\, dyds, \quad x\in\R^d, \, t>0,
\end{align}
where $p_t(x)=(2\pi t)^{-d/2}\exp(-|x|^2/2)$ is the transition probability
density of $B_t$.

Further, by Theorem \ref{th-3} (with $\sigma^2$ be replaced
by $(2\kappa)^{-1}\sigma^2$) and the relation (\ref{poisson-12}),
\begin{align}\label{th-16}
\log\E u_{-}(t,x)=\log\E u_{-}(t,0)
\sim\ \begin{cases}
\rho_1\Big({\theta\over 2\kappa}, {\sigma^2\over 2\kappa}\Big)t^{d/p}
& \hbox{in Regime I}\medskip \\
\rho_2\Big({\theta\over 2\kappa}, {\sigma^2\over 2\kappa}\Big)t^{3/2}
& \hbox{in Regime II}\medskip \\
\rho_3\Big({\theta\over 2\kappa}, {\sigma^2\over 2\kappa}\Big)
t^{4+d-2p\over 2}& \hbox{in Regime III}\medskip \\
\rho_4\Big({\theta\over 2\kappa}, {\sigma^2\over 2\kappa}\Big)t\log t& \hbox{in Regime IV}\medskip \\
\rho_5\Big({\theta\over 2\kappa}, {\sigma^2\over 2\kappa}\Big)t&
\hbox{in Regime V}\end{cases}\hskip.2in (t\to\infty).
\end{align}

By Theorem \ref{th-7}  and the relation (\ref{poisson-12}), under
$d/2<p<\min\{2, d\}$
\begin{align}\label{th-17}
\lim_{t\to\infty}{1\over t}\log\log \E u_+(t,x) ={2-p\over
4\kappa}\Big({p\over\sigma^2}\Big)^{p\over 2-p}
\Big(\theta\gamma(d,p)\Big)^{2\over 2-p}.
\end{align}
Under $p=2$ and $d=3$, we have
\begin{align}\label{th-18}
\lim_{t\to\infty}t^{-3/2}\log \E u_+(t,x)
=(2\kappa)^{3/2}\int_{\R^3}\E\Psi\bigg({\theta\over 2\kappa}
\int_0^1{ds\over \vert x+X_0(s)\vert^2}\bigg)dx
\end{align}
for every $\theta <\sigma^2/(16\kappa)$.

Our approach is guided by the idea known as Pascal principle,
which first introduced by Moreau et al (\cite{MOBC-1}, (\cite{MOBC-2})
and reformulated by Drewitz et al (\cite{DGRS}) in their context.
According to Pascal principle, the best way for the Brownian particle(s)
to avoid Poisson traps in the direct trapping reaction model or, to split
at the requested level
 in catalytic model, is to stay in a neighborhood of 0
(see Lemma \ref{moment-10} for some mathematical implementation of
 Pascal principle).
The mathematical formulation  is the relation
$$
\E\otimes\E_0\exp\bigg\{\pm \theta\int_0^t\ol{V}(s, B_s)ds\bigg\}\le \E\exp\bigg\{\pm \theta\int_0^t\ol{V}(s, 0)ds\bigg\}.
$$
 Different from \cite{MOBC-1}, (\cite{MOBC-2}
and \cite{DGRS},  our realization ( Lemma \ref{moment-10})
of Pascal principle relies on the Gaussian property (rather than
Markovian property) of the system.

Based on Pascal  principle and in
view of the equalities in Lemma \ref{moment-4},
much of our attention is on the asymptotic behaviors of the integrals
$$
\int_{\R^d}\psi\bigg(\theta\int_0^t{ds\over\vert x+X_0\vert^p}\bigg)dx
\hskip.1in \hbox{and}\hskip.1in
\int_{\R^d}\Psi\bigg(\theta\int_0^t{ds\over\vert x+X_0\vert^p}\bigg)dx.
$$
It should be pointed out that our treatment is very different in
the trapping and catalytic settings and in different regimes
marked by the combinations of $d$ and $p$.

The rest of this paper is organized as follows: In Section 3, we
establish moment identities and Pascal principle for random potential
function $\ol{V}(t, B_t)$. Both of these two results will play key roles in
the proofs of Theorems \ref{th-3} and \ref{th-7}
presented in Sections \ref{trap} and \ref{cata}, respectively. An appendix which
includes the proof of a technique result is given at the end of the
paper.

\section{Moment of occupation time with random potential}

As preparation for the proof of the main results, in this
section, we prove a key identity for the exponential moment of the
occupation time with random potential, and the Pascal principle. The
identity essentially says that the expectation (with respect to the
environment) of the exponential of the occupation time is the same
as the exponential of the expectation of a related random variable.
The Pascal principle essentially says that this expectation just mentioned in last sentence is
maximized when the particle does not move in the random medium.  This
principle is useful in deriving upper bound for the
asymptotic behavior of the occupation time.

According to Propositions 2.1 and 2.8 of \cite{CK} and
(\ref{poisson-8}), for any Borel measurable function $K(x)\ge 0$,
the compensated Poisson integral
\begin{align}\label{moment-2}
\ol{V}_K(t, x)=\int_{R^d}K(y-x)\big[\omega_t(dy)-dy\big]
\end{align}
is defined under the condition
\begin{align}\label{moment-3}
\int_{\R^d}\psi\Big(K(x)\Big)dx<\infty.
\end{align}

\begin{lemma}\label{moment-4} Let $K(x)\ge 0$ satisfy (\ref{moment-3}).
The time integral
\begin{align}\label{moment-5}
\int_0^t\ol{V}_K(s, B_s)ds
\end{align}
converges almost surely for any $t>0$ and
\begin{align}\label{moment-6}
&\E_0\otimes\E\exp\bigg\{-\int_0^t\ol{V}_K(s, B_s)ds\bigg\}\\
&=\E_0\exp\bigg\{\int_{\R^d}\E\psi\bigg(\int_0^tK(x+X_0(s)-B_s)ds\bigg)
dx\bigg\}.\nonumber
\end{align}

In addition, the equality
\begin{align}\label{moment-7}
&\E_0\otimes\E\exp\bigg\{\int_0^t\ol{V}_K(s, B_s)ds\bigg\}\\
&=\E_0\exp\bigg\{\int_{\R^d}\E\Psi\bigg(\int_0^tK(x+X_0(s)-B_s)ds\bigg)
dx\bigg\}\nonumber
\end{align}
holds.
\end{lemma}

\proof By (\ref{poisson-8}) and Fubini's theorem,
$$
\E_0\otimes\E\bigg(\int_0^t\vert \ol{V}_K(s,B_s)\vert ds\bigg)
=\E_0\int_0^t\E\vert \ol{V}_K(0,B_s)\vert ds=t\E\vert \ol{V}_K(0,
0)\vert<\infty,
$$
where the first step follows from (\ref{poisson-8}) and the last
step follows from Proposition 2.6 of \cite{CK}. So the integral in
(\ref{moment-5}) a.s. converges.

First, we claim that the
right hand side of  (\ref{moment-6}) is finite. Indeed, by Jensen inequality
\begin{eqnarray}\label{moment-8}
&&\int_{\R^d}\psi\bigg(\int_0^tK(x+X_0(s)-B_s)ds\bigg)
dx\\
&\le& {1\over t}\int_0^t\int_{\R^d}\psi\big(tK(x+X_0(s)-B_s)\big)dxds\nonumber\\
&=&\int_{\R^d}\psi\big(tK(x)\big)dx<\infty, \nonumber
\end{eqnarray}
where the last step follows from (\ref{moment-3}), monotonicity of $\psi(\cdot)$
and the fact that $\psi(2\theta)\le C\psi(\theta)$ ($\forall \theta>0$)
for some constant $C>0$.

To prove the equality in (\ref{moment-6}), we first consider the
case when $K(x)$ is bounded and locally supported. By
(\ref{poisson}) and Fubini theorem,
$$
\int_0^t\bigg[\int_{\R^d}K(x-B_s)\omega_s(dx)\bigg]ds
=\int_{\R^d}\bigg[\int_0^tK(X_y(s)-B_s)ds\bigg]\omega_0(dy)
$$
Hence, by the independence between $X_y(s)$ and $\omega_0(dy)$, and
by the independence among the Brownian motions $X_y(s)$, we have
\begin{eqnarray*}
&&\E\exp\bigg\{-\int_0^t\bigg[\int_{\R^d}K(x-B_s)\omega_s(dx)
\bigg]ds\bigg\}\\
&=&\E\exp\bigg\{\int_{\R^d}h(y)\omega_0(dy)\bigg\}\\
&=&\exp\bigg\{\int_{\R^d}\Big(e^{h(y)}-1\Big)dy\bigg\}\\
&=&\exp\bigg\{\int_{\R^d}
\bigg[\E\exp\bigg(-\int_0^tK(y+X_0(s)-B_s)ds\bigg)-1\bigg]dy\bigg\},
\end{eqnarray*}
where
$$
h(y)=\log\E\exp\bigg(-\int_0^tK(y+X_0(s)-B_s)ds\bigg).
$$
Consequently,
$$
\begin{aligned}
&\E\exp\bigg\{-\int_0^t\bigg[\int_{\R^d}K(x-B_s)\big[\omega_s(dx)-dx\big]
\bigg]ds\bigg\}\\
&=\exp\bigg\{\int_{\R^d}
\bigg[\E\exp\bigg(-\int_0^tK(y+X_0(s)-B_s)ds\bigg)-1+tK(y)\bigg]dy\bigg\}.
\end{aligned}
$$
Thus, (\ref{moment-6}) (with bounded and locally supported
$K(x)$) follows from the equality
\begin{align}\label{moment-1*}
\int_{\R^d}\bigg[\int_0^tK(y+X_0(s)-B_s)ds\bigg]dx=t\int_{\R^d}K(y)dy.
\end{align}

We now prove (\ref{moment-6}) with full generality. Let
$K_n(x)\ge 0$ be a non-decreasing sequence of bounded, locally supported
function such that $K_n(x)\uparrow K(x)$ point-wise. Notice that
$$
\E_0\otimes \E\bigg\vert\int_0^t \Big(\ol{V}_{K_n}(s,
B_s)-\ol{V}_{K}(s, B_s)\Big)ds\bigg\vert \le t\E_0\otimes
\E\Big\vert\ol{V}_{K_n}(0, 0)-\ol{V}_{K}(0, 0)\Big\vert.
$$
Using the decomposition in the proof of Proposition 2.6 of
\cite{CK}, one can show that the right hand side tends to zero as
$n\to\infty$.

We now take $n\to\infty$ on the both sides of
(\ref{moment-6}) with $K_n$ replacing $K$. The right hand side
passes from $K_n$ to $K$ by the monotonic convergence. On the other
hand, we note that
\begin{eqnarray*}
&&\E_0\otimes\E\exp\bigg\{-2\int_0^t\ol{V}_{K_n}(s, B_s)ds\bigg\}\\
&=&\E_0\exp\bigg\{\int_{\R^d}\E\psi\bigg(\int_0^t2K_n(x+X_0(s)-B_s)ds\bigg)
dx\bigg\}\\
&\le&\E_0\exp\bigg\{\int_{\R^d}\E\psi\bigg(\int_0^t2K(x+X_0(s)-B_s)ds\bigg)
dx\bigg\}\\
&<&\infty,
\end{eqnarray*}
where the last step follows from (\ref{moment-8}). Therefore, the
left hand side of (\ref{moment-6}) passes from $K_n$ to $K$ by the
uniform integrability of $\left\{\exp\left\{-\int_0^t\ol{V}_{K_n}(s,
B_s)ds\right\}:\; n\ge 1\right\}$. Hence, (\ref{moment-6}) holds for
$K$. The proof of (\ref{moment-7}) follows from a similar argument.
\qed

\begin{lemma}[Pascal principle]\label{moment-10}
Assume that $d/2<p<d$.
For any deterministic continuous function $b_s$ in $C\([0,t], \R^d\)$,
\begin{align}\label{moment-2*}
\int_{\R^d}\E\psi\bigg(\int_0^t{ds\over \vert
x+X_0(s)-b_s\vert^p}\bigg)dx\le\int_{\R^d}\E\psi\bigg(\int_0^t{ds\over \vert
x+X_0(s)\vert^p}\bigg)dx
\end{align}
and
\begin{align}\label{moment-11}
\E\int_{\R^d}\bigg[\int_0^t{ds\over \vert
x+X_0(s)-b_s\vert^p}\bigg]^mdx \le
\E\int_{\R^d}\bigg[\int_0^t{ds\over \vert x+X_0(s)\vert^p}\bigg]^m
dx,\hskip.2in m=2,3,\cdots.
\end{align}
Consequently (from (\ref{moment-11}), for any $\theta>0$
\begin{align}\label{moment-12}
\E\int_{\R^d}\Psi\bigg(\theta\int_0^t{ds\over \vert
x+X_0(s)-b_s\vert^p}\bigg)dx
\le\E\int_{\R^d}\Psi\bigg(\theta\int_0^t{ds\over \vert
x+X_0(s)\vert^p}\bigg)dx.
\end{align}
\end{lemma}

\proof Let $K(x)\ge 0$ be bounded, locally supported and radially
symmetric. As a corollary of Theorem 1.5, \cite{DSS},
$$
\begin{aligned}
&\int_{\R^d}
\bigg[\E\exp\bigg(-\int_0^tK(x+X_0(s)-b_s)ds\bigg)-1\bigg]dx\\
&\le \int_{\R^d}
\bigg[\E\exp\bigg(-\int_0^tK(x+X_0(s))ds\bigg)-1\bigg]dx.
\end{aligned}
$$
Subtracting
$$
t\int_{\R^d}K(x)dx
$$
on the both sides and by (\ref{moment-1*}) (with $B_s$ being replaced by
$b_S$ or 0),
$$
\int_{\R^d}\E\psi\bigg(\int_0^tK(x+X_0(s)-b_s)ds\bigg)dx
\le \int_{\R^d}
\E\psi\bigg(\int_0^tK(x+X_0(s))ds\bigg)dx.
$$
Taking a sequence of such function $K$'s approximating the function $|\cdot|^{-p}$, we obtain the desired inequality (\ref{moment-2*}).

Our approach for (\ref{moment-11}) relies on
Fourier transformation and the fact
that
\begin{align}\label{moment-13}
\int_{\R^d}{1\over\vert x\vert^p}e^{i\lambda\cdot x}dx =C{1\over
\vert \lambda\vert^{d-p}},
\end{align}
where $C>0$ is a constant.

Write
$$
h(y_1,\cdots, y_{m-1})=\int_{\R^d}{1\over\vert x\vert^p}
\prod_{k=1}^{m-1}{1\over\vert x+y_k\vert^p}dx.
$$
Then,
$$
\begin{aligned}
&\widehat{h}(\lambda_1,\cdots,\lambda_{m-1})
=\int_{(\R^d)^{m-1}}h(y_1,\cdots, y_{m-1})
\exp\Big\{i\sum_{k=1}^{m-1}\lambda_k\cdot y_k\Big\}dy_1\cdots dy_{m-1}\\
&=\int_{\R^d}{1\over\vert x\vert^p}
\exp\Big\{-i(\lambda_1+\cdots +\lambda_{m-1})\cdot x\Big\}dx
\prod_{k=1}^{m-1}\int_{\R^d}{1\over\vert y\vert^p}e^{i\lambda_k\cdot y}dy\\
&=C^m{1\over\vert \lambda_1+\cdots +\lambda_{m-1}\vert^{d-p}}
\prod_{k=1}^{m-1}{1\over\vert \lambda_k\vert^{d-p}}>0.
\end{aligned}
$$

Write
$$
\begin{aligned}
&\int_{\R^d}\bigg[\int_0^t{ds\over \vert x+X_0(s)-b_s\vert^p}\bigg]^mdx\\
&=\int_{[0,t]^m}\bigg[\int_{\R^d}
\prod_{k=1}^m{1\over \vert x+X_0(s_k)-b_{s_k}\vert^p}dx\bigg]ds_1\cdots ds_m\\
&=\int_{[0,t]^m}\bigg[\int_{\R^d}
{1\over\vert x\vert^p}\prod_{k=1}^m{1\over \vert x+\big(X_0(s_k)-X_0(s_m)\big)-
(b_{s_k}-b_{s_m})\vert^p}dx\bigg]ds_1\cdots ds_m\\
&=\int_{[0,t]^m}h\Big(\big(X_0(s_1)-X_0(s_m)\big)-
(b_{s_1}-b_{s_m}),\cdots,\\
&\hskip.8in  \big(X_0(s_{m-1})-X_0(s_m)\big)-
(b_{s_{m-1}}-b_{s_m})\Big)ds_1\cdots ds_m.
\end{aligned}
$$
By Fourier inversion,
$$
\begin{aligned}
&h\Big(\big(X_0(s_1)-X_0(s_m)\big)-
(b_{s_1}-b_{s_m}),\cdots,\big(X_0(s_{m-1})-X_0(s_m)\big)-
(b_{s_{m-1}}-b_{s_m})\Big)\\
&={1\over (2\pi)^{(m-1)d}}\int_{(\R^d)^{m-1}} d\lambda_1\cdots d\lambda_{m-1}
\widehat{h}(\lambda_1,\cdots,\lambda_{m-1})\\
&\times\exp\bigg\{
-i\sum_{k=1}^{m-1}\lambda_k\cdot\Big(\big(X_0(s_k)-X_0(s_m)\big)-
(b_{s_k}-b_{s_m})\Big)\bigg\}.
\end{aligned}
$$
Therefore
$$
\begin{aligned}
&\int_{\R^d}\bigg[\int_0^t{ds\over \vert x+X_0(s)-b_s\vert^p}\bigg]^mdx\\
&={1\over (2\pi)^{(m-1)d}}\int_{(\R^d)^{m-1}}d\lambda_1\cdots d\lambda_{m-1}
\widehat{h}(\lambda_1,\cdots,\lambda_{m-1})\int_{[0,t]^m}ds_1\cdots ds_m\\
&\times\exp\bigg\{
-i\sum_{k=1}^{m-1}\lambda_k\cdot\Big(\big(X_0(s_k)-X_0(s_m)\big)-
(b_{s_k}-b_{s_m})\Big)\bigg\}.
\end{aligned}
$$
Hence,
$$
\begin{aligned}
&\int_{\R^d}\E\bigg[\int_0^t{ds\over \vert x+X_0(s)-b_s\vert^p}\bigg]^mdx\\
&={1\over (2\pi)^{(m-1)d}}\int_{(\R^d)^{m-1}}d\lambda_1\cdots d\lambda_{m-1}
\widehat{h}(\lambda_1,\cdots,\lambda_{m-1})\int_{[0,t]^m}ds_1\cdots ds_m\\
&\times\exp\bigg\{-{1\over 2}
\Var\Big(\sum_{k=1}^{m-1}\lambda_k\cdot\Big(\big(X_0(s_k)-X_0(s_m)\big)\Big)
\bigg\}
\exp\bigg\{i\sum_{k=1}^{m-1}\lambda_k\cdot
(b_{s_k}-b_{s_m})\bigg\}\\
&\le{1\over (2\pi)^{(m-1)d}}\int_{(\R^d)^{m-1}}d\lambda_1\cdots d\lambda_{m-1}
\widehat{h}(\lambda_1,\cdots,\lambda_{m-1})\int_{[0,t]^m}ds_1\cdots ds_m\\
&\times\exp\bigg\{-{1\over 2}
\Var\Big(\sum_{k=1}^{m-1}\lambda_k\cdot\Big(\big(X_0(s_k)-X_0(s_m)\big)\Big)
\bigg\}\\
&=\int_{\R^d}\E\bigg[\int_0^t{ds\over \vert
x+X_0(s)\vert^p}\bigg]^mdx.
\end{aligned}
$$

The inequality (\ref{moment-12}) follows from the Taylor expansion
of $\Psi$ at $a=0$. \qed

\begin{remark}
The asymptotic version of Pascal principle (\ref{moment-2*}) for a different shape function has been obtained by \cite{PSSS}
and the full Pascal principle for that setting was obtained by \cite{PS}.
\end{remark}

Finally, we state a variety of the inequality (\ref{moment-11}) which will be used later.
For any $a>0$ write $K_a(x)=\vert x\vert^{-p}1_{\{\vert x\vert\ge a\}}$.
It is easy to check that
$$
\widehat{K_a}(\lambda)\equiv\int_{\R^d}K_a(x)e^{i\lambda\cdot
x}dx\ge 0 \hskip.2in\lambda\in\R^d.
$$
By the same argument as the one for (\ref{moment-11}),
\begin{align}\label{moment-14}
\int_{\R^d}\E\bigg[\int_0^tK_a\big(x+X_0(s)-b_s\big)ds\bigg]^mdx
\le \int_{\R^d}\E\bigg[\int_0^tK_a\big(x+X_0(s)\big)ds\bigg]^m
dx
\end{align}
for $m=2,3,\cdots$.

\section{Model of trapping reactions}\label{trap}

 We prove Theorem \ref{th-3} in this section. Throughout, we assume
that $d/2<p<d$. The integrability assertion (\ref{th-4}) follows
from (\ref{moment-6}) and (\ref{moment-8}) with $K(x)=\theta\vert
x\vert^{-p}$. In the following subsections, we establish the
asymptotics given in (\ref{th-5}) in different regimes.

The key is to estimate the exponent
\[\bar{\psi}_t(B)\equiv \int_{\R^d}\E
\psi\(\theta\int_0^t{ds\over\vert x+X_0(s)-B_s\vert^p}\)dx,\]
which depends on the initial obstacles (IO), moving obstacles (MO) and
the moving particle (MP).

It follows from Lemma \ref{moment-4} that
\begin{equation}\label{eq0731a}
\E\otimes\E_0\exp\bigg\{-\theta\int_0^tV(s, B_{s})ds\bigg\}=\E_0\exp\(\bar{\psi}_t(B)\).
\end{equation}
By Pascal principle, we have
\begin{equation}\label{eq0731b}
\log\E\otimes\E_0\exp\bigg\{-\theta\int_0^tV(s, B_{s})ds\bigg\}\le\bar{\psi}_t(0).
\end{equation}
So the proof of the upper bound is reduced to the estimate of the deterministic quantity $\bar{\psi}_t(0)$.

\subsection{Regime I: $p<2$}\label{R1}

 The main idea of the proof for this region is to consider the
asymptotic scales of the three components affecting $\bar{\psi}_t$.
From the shape function, we see that the IO is of the order
$t^{1/p}$ (see (\ref{eq0821a}) below). On the other hand, the MO and
MP are Brownian motions, and hence, of order $t^{1/2}$. When $p<2$,
 the MO and
MP are negligible.

By (\ref{moment-8}) (with $K(x)=\theta\vert x\vert^{-p}$),
\begin{equation}\label{eq0821a}
\bar{\psi}_t(0)\le\int_{\R^d}\psi\Big({t\theta\over\vert
x\vert^p}\Big)dx=(t\th)^{d/p}\int_{\R^d}\psi\Big({1\over\vert x\vert^p}\Big)dx.
\end{equation}
By Lemma
7.1 of \cite{Chen-1}, we have
\begin{align}\label{th-2}
\int_{\R^d}\psi\Big({1\over\vert x\vert^p}\Big)dx =\omega_d{p\over
d-p}\Gamma\Big({2p-d\over p}\Big),\hskip.3in \forall\ d/2<p<d.
\end{align}
By (\ref{eq0731b})-(\ref{th-2}), we obtain the upper
bound
$$
\limsup_{t\to\infty}t^{-d/p}
\log\E\otimes\E_0\exp\bigg\{-\theta\int_0^t\ol{V}(s, B_s)ds\bigg\} \le
\theta^{d/p}\omega_d{p\over d-p}\Gamma\Big({2p-d\over p}\Big).
$$

On the other hand, for $|x|>\de t^{\frac1p}$ and $|y|\le 2M\sqrt{t}$, we have as $t\to\infty$,
\[\frac{|x+y|}{|x|}\le1+\frac{|y|}{|x|}\le 1+\frac{2M}{\de}t^{\frac12-\frac1p}\to 1,\]
and hence, for $t$ large enough,
\[\frac{1}{|x+y|^p}\ge \frac{1-\de}{|x|^p}.\]
Therefore,
\begin{eqnarray*}
&&\E_0\exp\bigg\{\bar{\psi}_t(B)\bigg\}\\
&\ge&\E_0\int_{|x|\ge \de t^{1/p}}\E
\psi\(\theta\int_0^t{(1-\de)\over\vert x\vert^p}ds\)dx1_{\{\sup_{s\le t}(|B_s|,|X_0(s)|)\le M\sqrt{t}\}}\\
&\ge& \P_0\Big\{\max_{s\le t}\vert B_s\vert\le M\sqrt{t}\Big\}
\exp\bigg\{\P\Big\{\max_{s\le t}\vert X_0(s)\vert\le M\sqrt{t}\Big\}
\int_{\vert x\vert\ge\delta t^{1/p}}
\psi\bigg({\theta t (1-\delta)\over\vert x\vert^p}\bigg)dx\bigg\}\\
&=&\P_0\Big\{\max_{s\le t}\vert B_s\vert\le M\sqrt{t}\Big\}
\exp\bigg\{\big(\theta t (1-\delta)\big)^{d/p}
\P\Big\{\max_{s\le t}\vert X_0(s)\vert\le M\sqrt{t}\Big\}\\
&&\times\int_{\{\vert x\vert\ge\delta
(1-\delta)^{-1/p}\theta^{-1/p}\}} \psi\bigg({1\over\vert
x\vert^p}\bigg)dx\bigg\},
\end{eqnarray*}
for any $0<\delta<1$ and $M>0$, as $t$ is large.

The probabilities
$$
\P_0\Big\{\max_{s\le t}\vert B_s\vert\le M\sqrt{t}\Big\}\hskip.1in\hbox{and}
\hskip.1in\P\Big\{\max_{s\le t}\vert X_0(s)\vert\le M\sqrt{t}\Big\}
$$
can be made arbitrarily close to 1 by taking $M>0$ arbitrarily large. Thus,
$$
\begin{aligned}
&\liminf_{t\to\infty}t^{-d/p}\log\E_0\exp\bigg\{\bar{\psi}_t(B)\bigg\}\\
&\ge \big(\theta (1-\delta)\big)^{d/p}\int_{\{\vert x\vert\ge\delta
(1-\delta)^{-1/p}\theta^{-1/p}\}} \psi\bigg({1\over\vert
x\vert^p}\bigg)dx
\end{aligned}.
$$
Letting $\delta\to 0^+$ on the right hand side, the desired lower
bound follows from (\ref{moment-6}) and (\ref{th-2}). \qed

\subsection{Regime II: $p=2$ and $d=3$}\label{R2}

In this regime, the scales of the three components are the same so
the limit is obtained by a scaling argument.

By scaling, we have
\[\E_0\exp\bigg\{\bar{\psi}_t(B)\bigg\}
=\E_0\exp\bigg\{t^{3/2}\bar{\psi}_1(B)\bigg\}.
\]
By (\ref{eq0731a}) and (\ref{eq0731b}), all we need to show is the lower bound.

Note that the estimate
\begin{eqnarray*}
&&\E_0\exp\bigg\{t^{3/2}\bar{\psi}_1(B)\bigg\}\\
&\ge&\P_0\Big\{\max_{0\le s\le 1}\vert B_s\vert\le\epsilon\Big\}
\exp\(t^{3/2}\int_{\R^3}
\E\psi\bigg(\theta\int_0^1
{ds\over\big(|x+X_0(s)|+\epsilon\big)^2}\bigg)dx\),
\end{eqnarray*}
implies that
\[\liminf_{t\to\infty}t^{-3/2}\log\E_0\exp\bigg\{\bar{\psi}_t(B)\bigg\}
\ge \int_{\R^3}
\E\psi\bigg(\theta\int_0^1
{ds\over\big(|x+X_0(s)|+\epsilon\big)^2}\bigg)dx.\]
 Letting $\epsilon\to 0^+$ and applying Fatou's lemma we obtain the desired lower bound.
\qed

Next, we give an explicit upper bound for $\rho_2$.

\begin{lemma}\label{lem1001a}
\begin{align}\label{th-3a}
\rho_2(\theta,\sigma^2)\le {8\over 3}\pi^{3/2}\theta^{3/2}.
\end{align}
\end{lemma}
Proof:
The conclusion follows from the following estimate
$$
\begin{aligned}
&\int_{\R^3}
\psi\bigg(\theta\int_0^1{ds\over\vert x+X_0(s)\vert^2}\bigg)dx
\le\int_0^1\bigg[
\int_{\R^3}\psi\Big({\theta\over\vert x+X_0(s)\vert^2}\Big)dx\bigg]ds\\
&=\int_{\R^3}\psi\Big({\theta\over\vert x\vert^2}\Big)dx=\theta^{3/2}
\int_{\R^3}\psi\Big({1\over\vert x\vert^2}\Big)dx= {8\over 3}\pi^{3/2}\theta^{3/2},
\end{aligned}
$$
where the last step follows from (\ref{th-2}).
\qed

\subsection{Regime III: $2<p<{d+2\over 2}$}\label{R3}

In this regime, we consider the Taylor expansion of the function
$\psi$. It turns out that the quadratic term dominates the others.
The upper bound estimate is again obtained by (\ref{eq0731b}). The lower bound is obtained by
roughly making the particle motionless, namely, to restrict to the
event $\{\sup_{s\le t}|B_s|\le t^{1/p}\}$. Although the probability
of this event tends to 0, its exponential rate is less than that of
the scale determined in the upper bound.

By the fact that $\psi(a)\le 2^{-1} a^2$, ($a\ge 0$), we have
$$
\begin{aligned}
&\bar{\psi}_t(0) \le \int_{\R^d}{\theta^2\over
2}\E\bigg[\int_0^t{ds\over\vert x+X_0(s)\vert^p}
\bigg]^2dx\\
&={\theta^2\over 2}C(d,p)
\E\int_0^t\!\!\int_0^t{drds\over\vert X_0(r)-X_0(s)\vert^{2p-d}},
\end{aligned}
$$
where the second step follows from  (\ref{trap-1}).

Let $U$ be a $d$-dimensional standard normal random vector, i.e., $U\sim N(0, I_d)$. By Brownian scaling
$$
\begin{aligned}
&\E\int_0^t\!\!\int_0^t{drds\over\vert X_0(r)-X_0(s)\vert^{2p-d}}
=\sigma^{-(2p-d)}\E\vert U\vert^{-(2p-d)}\int_0^t\!\!\int_0^t
\vert r-s\vert^{-{2p-d\over 2}}drds\\
&=t^{4-2p+d\over 2}{2^{4+d-2p\over 2}d\omega_d\over
(2+d-2p)(4+d-2p)\sigma^{2p-d}}\pi^{-d/2} \Gamma(d-p).
\end{aligned}
$$

Summarizing our computation, we have obtained the desired upper bound
\begin{align}\label{trap-3}
\limsup_{t\to\infty}t^{-{4+d-2p\over 2}}\log
\E\otimes\E_0\exp\bigg\{-\theta\int_0^t\ol{V}(s, B_s)ds\bigg\}
\le\rho_3(\theta, \sigma^2).
\end{align}

On the other hand, for any $M>1$
\begin{align}\label{trap-4}
&\E_0\exp\bigg\{\bar{\psi}_t(B)\bigg\}\\
&\ge\E_0\exp\bigg\{\int_{\R^d}\E \psi\bigg(\theta\int_0^t{1_{\{\vert
x+X_0(s)\vert\ge Mt^{1/p}\}}
\over\vert x+X_0(s)-B_s\vert^p}ds\bigg)dx\bigg\}\nonumber\\
&\ge \P_0\Big\{\max_{s\le t}\vert B_s\vert\le t^{1/p}\Big\}
\exp\bigg\{\int_{\R^d}\E \psi\bigg({M^p\theta\over
(M+1)^p}\int_0^t{1_{\{\vert x+X_0(s)\vert\ge Mt^{1/p}\}} \over\vert
x+X_0(s)\vert^p}ds\bigg)dx\bigg\}.\nonumber
\end{align}

Given $\epsilon>0$ there is a $\delta>0$ such that
$\psi(a)\ge 2^{-1}(1-\epsilon)a^2$ whenever $0\le a\le\delta$.
Take $M>0$ sufficiently large so
$$
{M^p\theta\over (M+1)^p}\int_0^t{1_{\{\vert x+X_0(s)\vert\ge
Mt^{1/p}\}} \over\vert x+X_0(s)\vert^p}ds\le {M^p\theta\over
(M+1)^p}M^{-p}\le\delta.
$$
Consequently,
\begin{align}\label{trap-5}
&\int_{\R^d}\E
\psi\bigg({M^p\theta\over (M+1)^p}\int_0^t{1_{\{\vert x+X_0(s)\vert\ge Mt^{1/p}\}}
\over\vert x+X_0(s)\vert^p}ds\bigg)dx\\
&\ge {1\over 2}(1-\epsilon) \bigg({M^p\theta\over
(M+1)^p}\bigg)^2\int_{\R^d}\E\bigg[ \int_0^t{1_{\{\vert
x+X_0(s)\vert\ge Mt^{1/p}\}} \over\vert
x+X_0(s)\vert^p}ds\bigg]^2dx.\nonumber
\end{align}

Notice that
\begin{align}\label{trap-6}
&\int_{\R^d}\E\bigg[
\int_0^t{1_{\{\vert x+X_0(s)\vert\ge Mt^{1/p}\}}
\over\vert x+X_0(s)\vert^p}ds\bigg]^2dx\\
&=\int_0^t\!\!\int_0^t\E\bigg[\int_{\R^d}{1_{\{\vert x+X_0(r)\vert\ge Mt^{1/p}\}}
1_{\{\vert x+X_0(s)\vert\ge Mt^{1/p}\}}\over\vert x+X_0(r)\vert^p
\vert x+X_0(s)\vert^p}dx\bigg]drds\nonumber\\
&=\int_0^t\!\!\int_0^t\bigg[\int_{\{\vert x\vert\ge Mt^{1/p}\}}
{1\over\vert x\vert^p}\E{1_{\{\vert x+X_0(s)-X_0(r)\vert\ge Mt^{1/p}\}}\over
\vert x+X_0(s)-X_0(r)\vert^p}dx\bigg]drds.\nonumber
\end{align}

Recall that $U\sim N(0, I_d)$. By Brownian scaling, the right hand side of (\ref{trap-6})
becomes
$$
\begin{aligned}
&\int_0^t\!\!\int_0^t\bigg[\int_{\{\vert x\vert\ge Mt^{1/p}\}}
{1\over\vert x\vert^p}\E{1_{\{\vert x+\vert r-s\vert^{1/2}\sigma U\vert
\ge Mt^{1/p}\}}\over
\vert x+\vert r-s\vert^{1/2}\sigma U\vert^p}dx\bigg]drds\\
&=\sigma^{-(2p-d)}
\int_0^t\!\!\int_0^t\vert r-s\vert^{-{2p-d\over 2}}Q\Big(
M\sigma^{-1}t^{1/p}\vert r-s\vert^{-1/2}\Big)drds\\
&\ge\sigma^{-(2p-d)}Q\Big(M\sigma^{-1}u^{-1/2}t^{{1\over p}-{1\over
2}}\Big) \int\!\!\int_{\{[0,t]^2\cap\{\vert r-s\vert\ge ut\}\}}
\vert r-s\vert^{-{2p-d\over 2}}drds,
\end{aligned}
$$
where $Q(b)$ is defined by (\ref{trap-7}) in the Appendix.
Note that
$$
\int\!\!\int_{\{[0,t]^2\cap\{\vert r-s\vert\ge ut\}\}}
\vert r-s\vert^{-{2p-d\over 2}}drds
=t^{4-2p+d\over 2}\int\!\!\int_{\{[0,1]^2\cap\{\vert r-s\vert\ge u\}\}}
\vert r-s\vert^{-{2p-d\over 2}}drds.
$$
Therefore,
\begin{eqnarray}\label{eq0424a}
&&\liminf_{t\to\infty}t^{-{4-2p+d\over 2}}\int_{\R^d}\E\bigg[
\int_0^t{1_{\{\vert x+X_0(s)\vert\ge Mt^{1/p}\}}
\over\vert x+X_0(s)\vert^p}ds\bigg]^2dx\nonumber\\
&\ge&\sigma^{-(2p-d)}\lim_{b\to 0}Q(b)
\int\!\!\int_{\{[0,1]^2\cap\{\vert r-s\vert\ge u\}\}}
\vert r-s\vert^{-{2p-d\over 2}}drds.
\end{eqnarray}
Note that
\begin{equation}\label{eq0424b}
\lim_{u\to 0+}\int\!\!\int_{\{[0,1]^2\cap\{\vert r-s\vert\ge u\}\}}
\vert r-s\vert^{-{2p-d\over 2}}drds
={8\over
(2+d-2p)(4+d-2p)}.
\end{equation}
By (\ref{trap-5}), (\ref{eq0424a}), (\ref{eq0424b}) and Lemma \ref{lem0422a}, we get
$$
\liminf_{t\to\infty}t^{-{4-2p+d\over 2}} \int_{\R^d}\E
\psi\bigg({M^p\theta\over (M+1)^p}\int_0^t{1_{\{\vert
x+X_0(s)\vert\ge Mt^{1/p}\}} \over\vert x+X_0(s)\vert^p}ds\bigg)dx
\ge (1-\epsilon)\Big( {M\over M+1}\Big)^p\rho_3(\theta, \sigma^2).
$$
Bringing this back to (\ref{trap-4}) and by the classic estimate
$$
\log\P_0\Big\{\max_{s\le t}\vert B_s\vert\le t^{1/p}\Big\} \ge -C
t^{p-2\over p}=-o\Big(t^{4+d-2p\over 2}\Big)
$$
for large $t$ and taking  into account that $\epsilon$ can be arbitrarily small
and $M$ can be arbitrarily large, we obtain the desired lower bound.
\qed

\subsection{Regime V:
$p>\max\Big\{2, {d+2\over 2}\Big\}$}

The existence of the limit is obtained by super-additivity in $t$,
which implies the existence of the limit. To prove the
non-triviality of the limit, $\bar{\psi}$ is written into two part
by comparing the distance between obstacles and the particle with a
constant. When the distance is large, we have $\psi(u)\sim \frac12
u^2$. When the distance is small, we have $\psi(u)\sim u$. The upper bound can be estimated directly. The lower bound is
obtained using Jensen's inequality and the scaling property of the
Brownian motion.

It is straightforward to check that the function $\psi(\cdot)$ is
super-additive on $\R^+$: $\psi(\alpha
+\beta)\ge\psi(\alpha)+\psi(\beta)$ for $\alpha, \beta\ge 0$.
 In particular, for any
$t_1, t_2>0$,
$$
\begin{aligned}
&\int_{\R^d}\psi\bigg(\theta\int_0^{t_1+t_2}
{ds\over\vert x+X_0(s)-B_s\vert^p}\bigg)dx\\
&\ge \int_{\R^d}\psi\bigg(\theta\int_0^{t_1} {ds\over\vert
x+X_0(s)-B_s\vert^p}\bigg)dx
+\int_{\R^d}\psi\bigg(\theta\int_{t_1}^{t_1+t_2}{ds\over\vert
x+X_0(s)-B_s\vert^p} \bigg)dx.
\end{aligned}
$$
Notice that the two terms on the right hand side are independent
and,
$$
\int_{\R^d}\psi\bigg(\theta\int_{t_1}^{t_1+t_2}{ds\over\vert
x+X_0(s)-B_s\vert^p} \bigg)dx\buildrel d\over =
\int_{\R^d}\psi\bigg(\theta\int_0^{t_2} {ds\over\vert
x+X_0(s)-B_s\vert^p}\bigg)dx.
$$
So we have that
\[\E_0\exp\bigg\{\bar{\psi}_{t_1+t_2}(B)\bigg\}
\ge\E_0\exp\bigg\{\bar{\psi}_{t_1}(B)\bigg\}
\E_0\exp\bigg\{\bar{\psi}_{t_2}(B)\bigg\}.
\]
Therefore, the limit
$$
\lim_{t\to\infty}{1\over t}\log
\E_0\exp\bigg\{\bar{\psi}_{t}(B)\bigg\}
=\rho_5(\theta, \sigma^2)
$$
exists as an extended constant $0<\rho_5(\theta, \sigma^2)\le\infty$.

\begin{lemma}
There are two finite constants $C_1(\sigma^2)$ and
$C_2(\sigma^2)$ such that
\begin{equation}\label{eq0413a}
C_1(\sigma^2)\theta^{(d-2)/(p-2)}\le\rho_5(\theta, \sigma^2)\le
C_2(\sigma^2)\theta^{(d-2)/(p-2)}.\end{equation}
\end{lemma}
Proof:
Let $a>0$ and $0<\gamma <1$ be fixed but arbitrary. By Pascal principle and convexity,
\begin{eqnarray}\label{trap-9}
\bar{\psi}_{t}(B)\le\bar{\psi}_t(0)
&\le& \gamma\int_{\R^d}\E\psi\bigg(\theta\gamma^{-1}
\int_0^t{1\{\vert x+X_0(s)\vert\ge a\}
\over\vert x+X_0(s)\vert^p}ds\bigg)dx\nonumber\\
&&+(1-\gamma) \int_{\R^d}\E\psi\bigg({\theta\over 1-\gamma}
\int_0^t{1\{\vert x+X_0(s)\vert <a\} \over\vert
x+X_0(s)\vert^p}ds\bigg)dx.
\end{eqnarray}

For the first term on the right hand side
\begin{align}\label{trap-10}
&\int_{\R^d}\E\psi\bigg(\theta\gamma^{-1}
\int_0^t{1\{\vert x+X_0(s)\vert\ge a\}
\over\vert x+X_0(s)\vert^p}ds\bigg)dx\\
&\le {\theta^2\over 2\gamma^2}\int_{\R^d}\E\bigg[
\int_0^t{1\{\vert x+X_0(s)\vert\ge a\}
\over\vert x+X_0(s)\vert^p}ds\bigg]^2dx.\nonumber
\end{align}
By  the computation next to (\ref{trap-6}),
\begin{align}\label{trap-10'}
&\int_{\R^d}\E\bigg[
\int_0^t{1\{\vert x+X_0(s)\vert\ge a\}
\over\vert x+X_0(s)\vert^p}ds\bigg]^2dx\\
&=\int_0^t\!\!\int_0^t\bigg[\int_{\{\vert x\vert\ge a\}}
{1\over\vert x\vert^p}\E{1_{\{\vert x+X_0(s)-X_0(r)\vert\ge a\}}\over
\vert x+X_0(s)-X_0(r)\vert^p}dx\bigg]drds\nonumber\\
&=\sigma^{-(2p-d)} \int_0^t\!\!\int_0^t
\vert r-s\vert^{-{2p-d\over 2}} Q\Big(a\sigma^{-1}\vert
r-s\vert^{-1/2}\Big)drds,\nonumber
\end{align}
where $Q(b)$ is defined in
(\ref{trap-7}), and the equation
follows from Brownian scaling and variable substitution.

Notice that
\begin{align}\label{trap-13a}
&\int_0^t\!\!\int_0^t \vert r-s\vert^{-{2p-d\over 2}}
Q\Big(a\sigma^{-1}\vert r-s\vert^{-1/2}\Big)drds =2\int_0^t
(t-s)s^{-{2p-d\over 2}}Q\Big(a\sigma^{-1}s^{-1/2}\Big)ds.
\end{align}
According to Lemma \ref{a-1}, $Q(b)=O(b^{-(2p-d)})$ as $b\to\infty$. By
(\ref{trap-8})
$$
\int_0^\infty s^{-{2p-d\over 2}}Q\Big(a\sigma^{-1}s^{-1/2}\Big)ds<\infty
$$
as $p>{d+2\over 2}$. Consequently, as $t\to\infty$,
$$
\begin{aligned}
&\int_0^t (t-s)s^{-{2p-d\over 2}}Q\Big(a\sigma^{-1}s^{-1/2}\Big)ds
\sim t\int_0^\infty s^{-{2p-d\over 2}}Q\Big(a\sigma^{-1}s^{-1/2}\Big)ds\\
&=2\Big({\sigma\over a}\Big)^{2p-d-2}t\int_0^\infty
s^{2p-d-3}Q(s)ds.
\end{aligned}
$$

Summarizing the computation, we have
\begin{align}\label{trap-11}
&\int_{\R^d}\E\bigg[
\int_0^t{1\{\vert x+X_0(s)\vert\ge a\}
\over\vert x+X_0(s)\vert^p}ds\bigg]^2dx\\
&\le \big(1+o(1)\big)a^{-(2p-d-2)}D(d,p,\sigma^2) t \nonumber
\end{align}
as $t\to\infty$, where
$$
D(d,p,\sigma^2)=4\sigma^{-2}\int_0^\infty s^{2p-d-3}Q(s)ds.
$$

As for the second term on the right hand side of (\ref{trap-9}),
taking $K(x)=\theta (1-\gamma)^{-1}\vert x\vert^{-p}1_{\{\vert x\vert <a\}}$ in
(\ref{moment-8}) gives
\begin{align}\label{trap-12}
&\int_{\R^d}\E\psi\bigg({\theta\over 1-\gamma}
\int_0^t{1\{\vert x+X_0(s)\vert <a\}
\over\vert x+X_0(s)\vert^p}ds\bigg)dx\\
&\le \int_{\{\vert x\vert<a\}}
\psi\bigg({\theta\over 1-\gamma}{t\over\vert x\vert^p}\bigg)dx\nonumber\\
&=\Big({\theta t\over 1-\gamma}\Big)^{d/p}
\int_{\{\vert x\vert\le a(1-\gamma)^{1/p}(\theta t)^{-1/p}\}}
\psi\Big({1\over\vert x\vert^p}\Big)dx\nonumber\\
&\sim\Big({\theta t\over 1-\gamma}\Big)^{d/p}
\int_{\{\vert x\vert\le a(1-\gamma)^{1/p}(\theta t)^{-1/p}\}}{1\over\vert x\vert^p}dx
\nonumber\\
&=\Big({\theta t\over 1-\gamma}\Big)^{d/p}{d\omega_d\over d-p}a^{d-p}
\Big((1-\gamma)^{1/p}(\theta t)^{-1/p}\Big)^{d-p}\nonumber\\
&={\theta t\over 1-\gamma}{d\omega_d\over d-p}a^{d-p},\hskip.2in
(t\to\infty). \nonumber
\end{align}

Combining (\ref{trap-10}), (\ref{trap-11}) and (\ref{trap-12}), we
have
\begin{align}\label{trap-13}
&\limsup_{t\to\infty}{1\over t}\log
\E_0\exp\bigg\{\bar{\psi}_{t}(B)\bigg\}\\
&\le (2\gamma)^{-1}\theta^2 a^{-(2p-d-2)}D(d,p,\sigma^2)+{\theta
d\omega_d\over d-p}a^{d-p}.\nonumber
\end{align}

Taking $\gamma\to 1^-$ and choosing $a$ to minimize the right hand
side of (\ref{trap-13}), we see that there is a constant
$C_2(\sigma^2)$ (drop its dependency on $d,\;p$) such that
$$
\rho_5(\theta, \sigma^2)\le
C_2(\sigma^2)\theta^{(d-2)/(p-2)}.
$$
 This proves the upper bound in
(\ref{eq0413a}), and hence, the finiteness of $\rho_5(\theta,
\sigma^2)$. In fact, by elementary calculus, it is easy to show that
$C_2(\sigma^2)=K(d,p)\sigma^{-2(d-p)/(p-2)}$ with
\[K(d,p)=\left\{\(\frac{d-p}{2p-d-2}\)^{\frac{2p-d-2}{p-2}}+\(\frac{2p-d-2}{p-2}\)^{\frac{2p-d-2}{d-p}}\right\}
D_1(d,p)^{\frac{2p-d-2}{d-p}}\(\frac{d\omega_d}{d-p}\)^{\frac{2p-d-2}{p-2}},\]
where $D_1(d,p)=\frac12 \sigma^2 D(d,p,\sigma^2)$.

On the other hand, let $\alpha$ be a constant to be decided later.
By Brownian scaling, we have
\begin{eqnarray*}
&&\frac{1}{t}\log\E_0\exp\left\{\bar{\psi}_{t}(B)\right\}\\
&=&\frac{1}{t}\E_0\int_{\R^d}\E\psi\bigg(\theta\int_0^t{ds\over
\vert
x+X_0(s)-B_s\vert^p}\bigg)dx\\
&=&\frac{1}{t}\E_0\otimes\E\int_{\R^d}\psi\bigg(\theta^{1-\alpha}\int_0^{\theta^\alpha
t}{ds\over
\vert x+\theta^{-\alpha/2}(X_0(s)-B_s)\vert^p}\bigg)dx\\
&=&\frac{1}{t}\E_0\otimes\E\int_{\R^d}\psi\bigg(\theta^{1-\alpha}\theta^{\alpha
p/2}\int_0^{\theta^\alpha t}{ds\over
\vert x\theta^{\alpha/2}+X_0(s)-B_s\vert^p}\bigg)dx\\
&=&\theta^{-\alpha d/2}\theta^\alpha \frac{1}{\theta^\alpha
t}\E_0\otimes\E\int_{\R^d}\psi\bigg(\theta^{1-\alpha}\theta^{\alpha
p/2}\int_0^{\theta^\alpha t}{ds\over \vert
x+X_0(s)-B_s\vert^p}\bigg)dx.
\end{eqnarray*}
Taking $\alpha=-\frac{2}{p-2}$ and letting $t\to\infty$, we see that
$$
\begin{aligned}
\rho_5(\theta, \sigma^2)&\ge
\theta^{(d-2)/(p-2)}\lim_{t\to\infty}\frac{1}{t}
\E_0\otimes\E\int_{\R^d}\psi\bigg(\int_0^{t}{ds\over
\vert x+X_0(s)-B_s\vert^p}\bigg)dx\\
&\ge\theta^{(d-2)/(p-2)}C_1(\sigma^2)
\end{aligned}
$$
where the existence of the limit follows from the superadditivity,
and
\[C_1(\sigma^2)=\E_0\otimes\E\int_{\R^d}\psi\bigg(\int_0^{1}{ds\over
\vert x+X_0(s)-B_s\vert^p}\bigg)dx.\]
\qed

\subsection{Regime IV: $p={d+2\over 2}>2$}

One of the differences of Regimes IV and V is that the second term on the right hand side of (\ref{eq0424c}) does make contribution in Regime IV.  A much more delicate
treatment is needed for this critical case. We first prove the upper bound.

\begin{proposition}\label{prop0422a}
\begin{align}\label{trap-13'}
\limsup_{t\to\infty}{1\over t\log t}\log \E_0\exp\bigg\{
\bar{\psi}_{t}(B)\bigg\} \le\rho_4(\theta, \sigma^2).
\end{align}
\end{proposition}
Proof: We continuing the calculation of (\ref{trap-13a}) of the
last subsection
Setting $p={d+2\over 2}$, we consider the decomposition
$$
\begin{aligned}
&\int_0^t (t-s)s^{-1}Q\Big(a\sigma^{-1}s^{-1/2}\Big)ds\\
&=\int_0^b (t-s)s^{-1}Q\Big(a\sigma^{-1}s^{-1/2}\Big)ds
+\int_b^t (t-s)s^{-1}Q\Big(a\sigma^{-1}s^{-1/2}\Big)ds\\
&=O(t)+\Big(1+
q(b)\Big)C(d,p)2^{d-p-1}d\omega_d(2\pi)^{-d/2}\Gamma(d-p)
\int_b^t(t-s)s^{-1}ds\\
&=O(t)+\Big(1+
q(b)\Big)C(d,p)2^{d-p-1}d\omega_d(2\pi)^{-d/2}\Gamma(d-p) t\log
t,\hskip.2in (t\to\infty),
\end{aligned}
$$
where $q(b)\to 0$ as $b\to \infty$ and the second step follows from
(\ref{trap-8}).

Together with (\ref{trap-10'}) and (\ref{trap-13a})
in last subsection, we have,
\begin{align}\label{trap-11a}
&\int_{\R^d}\E\bigg[ \int_0^t{1\{\vert x+X_0(s)\vert\ge a\}
\over\vert x+X_0(s)\vert^p}ds\bigg]^2dx\\
&=\big(1+o(1)\big)2\sigma^{-2}
C(d,p)2^{d-p-1}d\omega_d(2\pi)^{-d/2}\Gamma(d-p) t\log t ,\nonumber
\end{align}
as $t\to\infty$.

The computation of the constant needs a little attention in this regime.
We  claim that
\begin{align}\label{trap-12a}
\theta^2\sigma^{-2}
C(d,p)2^{d-p-1}d\omega_d(2\pi)^{-d/2}\Gamma(d-p)=\rho_4(\theta,\sigma^2)
\end{align}

By (\ref{trap-2}) with the relation $p={d+2\over 2}$, the left hand
side of (\ref{trap-12a}) is equal to
$$
2^{d-4\over 2}d\omega_d\Big({\theta\over\sigma}\Big)^2
\bigg({\displaystyle
\Gamma\Big({d-p\over 2}\Big)\over\displaystyle
\Gamma\Big({p\over 2}\Big)}\bigg)^2
=2^{d-4\over 2}d\omega_d\Big({\theta\over\sigma}\Big)^2
\bigg({\displaystyle
\Gamma\Big({d-2\over 4}\Big)\over\displaystyle
\Gamma\Big({d-2\over 4}+1\Big)}\bigg)^2
$$
Therefore, our assertion follows from the relation that
$$
\Gamma\Big({d-2\over 4}+1\Big)={d-2\over 4}\Gamma\Big({d-2\over 4}\Big).
$$

In view of (\ref{trap-10}), (\ref{trap-12}) and (\ref{trap-11a}),
\[
\limsup_{t\to\infty}{1\over t\log t}\log \E_0\exp\bigg\{
\bar{\psi}_{t}(B)\bigg\} \le\gamma^{-1}\rho_4(\theta,
\sigma^2).
\]
Letting $\gamma\to 1^-$ on the right hand side, we finish the proof of (\ref{trap-13'}).
\qed

It remains to establish the corresponding lower bound. Again, the
truncation level $a>0$ is fixed but arbitrary. The challenge is to
reverse the inequality by Taylor expansion given in (\ref{trap-10}).
Write $K_a(x)=\vert x\vert^{-p}1_{\{\vert x\vert\ge a\}}$ and
$$
v(t,x)=\E\exp\bigg\{-\theta\int_0^tK_a\big(x+X_0(s)\big)ds\bigg\},
\hskip.2in (t,x)\in\R^+\times\R^d.
$$

\begin{lemma}
\begin{align}\label{trap-16}
\lim_{t\to\infty}{1\over t\log t}
\int_{\R^d}\E\psi\bigg(\theta\int_0^tK_a\big(x+X_0(s)\big)ds\bigg)dx
=\rho_4(\theta, \sigma^2).
\end{align}
\end{lemma}
Proof:
By the fundamental theorem of calculus, we have
$$
\begin{aligned}
&\exp\bigg\{-\theta\int_0^tK_a\big(x+X_0(s)\big)ds\bigg\}\\
&=1-\theta\int_0^tK_a\big(x+X_0(s)\big)\exp\bigg\{-\theta\int_s^t
K_a\big(x+X_0(u)\big)du\bigg\}ds.
\end{aligned}
$$
Taking expectation on the both sides, by Markov property, we have
\begin{align}\label{trap-14}
v(t,x)=1-\theta\sigma^{-d}\int_0^t\int_{\R^d}p_{t-s}
\Big({x-y\over\sigma}\Big)K_a(y)v(s,y)dyds,
\end{align}
where
$$
p_t(x)={1\over (2\pi t)^{d/2}}\exp\Big\{-{\vert x\vert^2\over
2t}\Big\}.
$$
Recall that $\psi(x)=e^{-x}-1+x$. Hence,
$$
\begin{aligned}
\E\psi\bigg(\theta\int_0^tK_a\big(x+X_0(s)\big)ds\bigg)&
=v(t, x)-1+\theta\E\int_0^tK_a\big(x+X_0(s)\big)ds\\
&=\theta\sigma^{-d}
\int_{\R^d}K_a(y)\int_0^tp_{t-s}\Big({x-y\over\sigma}\Big)\Big(1-v(s,y)\Big)dyds.
\end{aligned}
$$
Therefore
$$
\begin{aligned}
&\int_{\R^d}\E\psi\bigg(\theta\int_0^tK_a\big(x+X_0(s)\big)ds\bigg)dx\\
&=\theta\int_{\R^d}K_a(y)\int_0^t\Big(1-v(s,y)\Big)dyds\\
&=\theta^2\sigma^{-d}\int_{\R^d}K_a(y)\int_0^t\bigg[\int_0^s
\int_{\R^d}p_{s-u}\Big({y-x\over\sigma}\Big)K_a(x)v(u,x)dxdu
\bigg]dyds\\
&=\theta^2\sigma^{-d}
\int_{\R^d}\int_{\R^d}K_a(x)K_a(y)\bigg[\int_0^tv(u,x)\bigg(\int_0^{t-u}
p_s\Big({y-x\over\sigma}\Big)ds\bigg)du\bigg]dxdy,
\end{aligned}
$$
where the second step follows from (\ref{trap-14}).

Taking Laplace transform on the both sides, for any $\lambda>0$
$$
\begin{aligned}
&\int_0^\infty e^{-\lambda t}\bigg[
\int_{\R^d}\E\psi\bigg(\theta\int_0^tK_a\big(x+X_0(s)\big)ds\bigg)dx\bigg]dt\\
&=\lambda^{-1}\theta^2\sigma^{-d}\int\!\!\int_{\R^d\times\R^d}K_a(x)K_a(y)
\bigg(\int_0^\infty e^{-\lambda t}v(t,x)dt\bigg)\bigg(\int_0^\infty
e^{-\lambda t} p_t\Big({y-x\over\sigma}\Big)dt\bigg)dxdy.
\end{aligned}
$$

By the inequality $e^{-x}\ge 1-x$, ($x\ge 0$),
$$
v(t,x)\ge 1-\theta\E\int_0^tK_a(x+X_0(s))ds
=1-\theta\sigma^{-d}\int_{\R^d}\int_0^tK_a(z)p_s\Big({z-x\over\sigma}\Big)dsdz
$$
Hence,
$$
\begin{aligned}
\int_0^\infty e^{-\lambda t}v(t,x)dt&\ge\lambda^{-1}\bigg\{1-\theta
\sigma^{-d} \int_{\R^d}K_a(z)\bigg(\int_0^\infty e^{-\lambda t}p_t
\Big({z-x\over\sigma}\Big)dt\bigg)dz\bigg\}\\
&\ge\lambda^{-1}\bigg\{1-\theta \sigma^{-d}
\int_{\R^d}K_a(z)\bigg(\int_0^\infty p_t
\Big({z-x\over\sigma}\Big)dt\bigg)dz\bigg\}.
\end{aligned}
$$

Notice that $d\ge 3$. It is well-known that
$$
\int_0^\infty p_t(x)dt=C_d\vert x\vert^{-(d-2)}\hskip.2in x\in\R^d.
$$
Summarizing our estimate,
\begin{align}\label{trap-15}
&\int_0^\infty e^{-\lambda t}\bigg[
\int_{\R^d}\E\psi\bigg(\theta\int_0^tK_a\big(x+X_0(s)\big)ds\bigg)dx\bigg]dt\\
&\ge\lambda^{-2}\theta^2\sigma^{-d}
\bigg\{\int_{\R^d}\int_{\R^d}K_a(x)K_a(y)
\bigg(\int_0^\infty e^{-\lambda t}
p_t\Big({y-x\over\sigma}\Big)dt\bigg)dxdy\nonumber\\
&-C\int_{\R^d}\int_{\R^d}\int_{\R^d}K_a(x)K_a(y)K_a(z){1\over\vert
y-x\vert^{d-2}} {1\over\vert z-x\vert^{d-2}}dxdydz\bigg\}.\nonumber
\end{align}

For the first term on the right hand side, notice that
$$
\begin{aligned}
&\int\!\!\int_{\R^d\times\R^d}K_a(x)K_a(y)
\bigg(\int_0^\infty e^{-\lambda t}
p_t\Big({y-x\over\sigma}\Big)dt\bigg)dxdy\\
&=\sigma^{2d}\int_0^\infty e^{-\lambda t}t^{d/2}\bigg[\int\!\!\int_{\R^d\times\R^d}
K_a(\sigma t^{1/2}x)K_a(\sigma t^{1/2}y)p_1(y-x)
dxdy\bigg]dt\\
&=\sigma^{2(d-p)}\int_0^\infty e^{-\lambda t} t^{-1} Q(a\sigma^{-1} t^{-1/2})dt\\
&\sim \sigma^{2(d-p)}C(d,p)2^{d-p-1}d\omega_d(2\pi)^{-d/2}\Gamma(d-p)
\log{1\over\lambda}
\end{aligned}
$$
as $\lambda\to 0^+$, where $Q(b)$ is defined in (\ref{trap-7}) and the last
step follows from (\ref{trap-8}) and the fact (Lemma \ref{a-1}) that
$t^{-1} Q(a\sigma t^{-1/2})$ is bounded for $t$ in a neighborhood of 0.

As for the second term in (\ref{trap-15}), notice that
$$
\int_{\R^d}\int_{\R^d}K_a(y){1\over \vert y-x\vert^{d-2}}dy
\le\int_{\R^d}\!\int_{\R^d}{1\over \vert y\vert^p\vert
y-x\vert^{d-2}}dy =C_1{1\over \vert x\vert^{p-2}},\hskip.2in
x\in\R^d.
$$
Therefore,
$$
\begin{aligned}
&\int_{\R^d}\int_{\R^d}\int_{\R^d}
K_a(x)K_a(y)K_a(z){1\over\vert y-x\vert^{d-2}}
{1\over\vert z-x\vert^{d-2}}dxdydz\\
&\le C_2\int_{\R^d}K_a(x){1\over\vert x\vert^{p-2}}{1\over\vert
x\vert^{p-2}}dx = C_2\int_{\vert x\vert\ge a\}}{dx\over\vert
x\vert^{3p-4}}<\infty,
\end{aligned}
$$
where the last step follows from the fact that $3p-4>d$ in Regime IV.

By (\ref{trap-15}), with $p={d+2\over 2}$
$$
\begin{aligned}
&\liminf_{\lambda\to 0^+}\lambda^2\Big(\log{1\over\lambda}\Big)^{-1}
\int_0^\infty e^{-\lambda t}\bigg[
\int_{\R^d}\E\psi\bigg(\theta\int_0^tK_a\big(x+X_0(s)\big)ds\bigg)dx\bigg]dt\\
&\ge \theta^2\sigma^{-2}C(d,p)2^{d-p-1}d\omega_d(2\pi)^{-d/2}\Gamma(d-p)
=\rho_4(\theta, \sigma^2).
\end{aligned}
$$
where the equality comes from (\ref{trap-12a}).

On the other hand, from the relation
$$
\int_{\R^d}\E\psi\bigg(\theta\int_0^tK_a\big(x+X_0(s)\big)ds\bigg)dx
\le {\theta^2\over
2}\int_{\R^d}\E\bigg[\int_0^tK_a\big(x+X_0(s)\big)ds\bigg]^2dx
$$
and from  (\ref{trap-11})
 we derive that
$$
\limsup_{\lambda\to 0^+}\lambda^2\Big(\log{1\over\lambda}\Big)^{-1}
\int_0^\infty e^{-\lambda t}\bigg[
\int_{\R^d}\E\psi\bigg(\theta\int_0^tK_a\big(x+X_0(s)\big)ds\bigg)dx\bigg]dt
\le\rho_4(\theta, \sigma^2).
$$
Consequently,
$$
\lim_{\lambda\to 0^+}\lambda^2\Big(\log{1\over\lambda}\Big)^{-1}
\int_0^\infty e^{-\lambda t}\bigg[
\int_{\R^d}\E\psi\bigg(\theta\int_0^tK_a\big(x+X_0(s)\big)ds\bigg)dx\bigg]dt
=\rho_4(\theta, \sigma^2).
$$
(\ref{trap-16}) then follows from Tauberian theorem (see Lemma 2.1.1. of Yakimiv \cite{yakimiv},
p91).
\qed

Finally, we are ready to consider the lower bounded in this regime.

\begin{proposition}\label{prop0422b}
\begin{equation}\label{eq0422a}
\lim_{t\to\infty}{1\over t\log t}\log \E_0\exp\bigg\{
\bar{\psi}_{t}(B)\bigg\} =\rho_4(\theta, \sigma^2).
\end{equation}
\end{proposition}
Proof:
Let $a>1$. We have
$$
\begin{aligned}
&\E_0\exp\(
\bar{\psi}_{t}(B)\)\\
&\ge\E_0\exp\bigg\{\int_{\R^d}\E\psi\bigg(\theta\int_0^t
K_a\big(x+X_0(s)-B_s\big)ds\bigg)dx\bigg)\bigg\}\\
&\ge\P_0\Big\{\max_{s\le t}\vert B_s\vert\le 1\Big\}
\exp\bigg\{\int_{\R^d}\E\psi\bigg(\Big({a\over
a+1}\Big)^p\theta\int_0^t K_a\big(x+X_0(s)\big)ds\bigg)dx\bigg\}.
\end{aligned}
$$
Replacing $\theta$ by $\Big({a\over a+1}\Big)^p\theta$ in (\ref{trap-16}),
$$
\liminf_{t\to\infty}{1\over t\log t}\log \E_0\exp\bigg\{
\bar{\psi}_{t}(B) \bigg\}\ge \rho_4\Big(\Big({a\over
a+1}\Big)^p\theta, \sigma^2\Big).
$$
Letting $a\to\infty$ on the right hand side and then combining it with
(\ref{trap-13'}), we obtain (\ref{eq0422a}).
\qed

Finally, the conclusion of Theorem \ref{th-3} for Regime IV follows from (\ref{moment-6}) with
$K(x)=\theta\vert x\vert^{-p}$ and Propositions \ref{prop0422a} and \ref{prop0422b}.

\section{Brownian motion in catalytic medium}\label{cata}

We prove Theorem \ref{th-7} in this section. Note that the
Assumption (\ref{th-1}) remains in force in this section. The proof
is splitted into three sub-sections according to the value of $p$.

The following notations are used in this section. For $R>0$ and
$x\in\R^d$, $B(x,R)$ represents the $d$-dimensional ball with the
center $x$ and the radius $R$. Given an open domain $D\subset \R^d$,
$W^{1,2}(D)$ is the Sobolev space over $D$, defined as the closure
of the inner product space of the infinitely differentiable
functions compactly supported in $D$ under the Sobolev norm
$$
\|g\|_H=\Big\{\|g\|_{{\cal L}^2(D)}^2+\|\nabla g\|_{{\cal L}^2(D)}^2\Big\}^{1/2}
$$
Write
$$
{\cal F}_d(D)=\Big\{g\in W^{1,2}(D);\hskip.1in \|g\|_{{\cal L}^2(D)}=1\Big\}
$$
In particular, ${\cal F}_d={\cal F}_d(\R^d)$.

\subsection{Sub-critical case $p<2$}\label{c-1}

 The strategy of the proof is as follows: For the upper bound, a
direct use of Pascal principle enables us to assume the particle is
motionless. For the lower bound, we restrict ourselves to the event
that the particle is almost motionless in the sense of
$\left\{\sup_{s\le t}|B_s|\le\epsilon\right\}$. Although the
probability of this event tends to 0 as $t\to\infty$, its rate is
much slower than double exponential.

Given $0<\gamma <1$, by convexity of $\Psi(\cdot)$, we have
\begin{align}\label{cata-1}
&\int_{\R^d}
\E\Psi\bigg(\theta\int_0^t{ds\over \vert x+X_0(s)\vert^p}\bigg)dx\\
&\le (1-\gamma) \int_{\R^d}\E\Psi\bigg({\theta\over 1-\gamma}
\int_0^t{1\{\vert x+X_0(s)\vert> t^{1/p}\}\over \vert x+X_0(s)\vert^p}ds\bigg)dx
\nonumber\\
&+\gamma \int_{\R^d}\E\Psi\bigg(\gamma^{-1}\theta \int_0^t{1\{\vert
x+X_0(s)\vert\le t^{1/p}\}\over \vert x+X_0(s)\vert^p}ds
\bigg)dx.\nonumber
\end{align}

For the first term, by Jensen inequality, we get
\begin{eqnarray}\label{eq0422b}
&&\int_{\R^d}\Psi\bigg({\theta\over 1-\gamma}
\int_0^t{1\{\vert x+X_0(s)\vert> t^{1/p}\}\over \vert x+X_0(s)\vert^p}ds\bigg)dx\nonumber\\
&\le& {1\over t}\int_0^t\int_{\R^d}\Psi\bigg({\theta\over 1-\gamma}
{t1\{\vert x+X_0(s)\vert> t^{1/p}\}\over \vert x+X_0(s)\vert^p}ds\bigg)dxds\nonumber\\
&=&\int_{\{\vert x\vert\ge t^{1/p}\}}\Psi\bigg({\theta\over 1-\gamma}
{t\over \vert x\vert^p}\bigg)dx\nonumber\\
&=&\Big({\theta t\over 1-\gamma}\Big)^{d/p} \int_{\{\vert x\vert\ge
(1-\gamma)^{1/p}\theta^{-1/p}\}}\Psi\bigg( {1\over \vert
x\vert^p}\bigg)dx<\infty.
\end{eqnarray}

As for the second term of (\ref{cata-1}), we use the bound
$$
\int_0^t{1\{\vert x+X_0(s)\vert\le t^{1/p}\}\over \vert
x+X_0(s)\vert^p} ds\le 1_{\{\vert x\vert\le t^{1/p}+\max_{s\le
t}\vert X_0(s)\vert\}} \sup_{x\in\R^d}\int_0^t{ds\over \vert
x+X_0(s)\vert^p}.
$$
Consequently,
\begin{eqnarray}\label{eq0422c}
&&\int_{\R^d}\E\Psi\bigg(\gamma^{-1}\theta
\int_0^t{1\{\vert x+X_0(s)\vert\le t^{1/p}\}\over \vert x+X_0(s)\vert^p}
ds\bigg)dx\\
&\le& \omega_d\E\Bigg[\Big(t^{1/p}+\max_{s\le t}\vert X_0(s)\vert\Big)^d
\Psi\bigg(\gamma^{-1}\theta\sup_{y\in\R^d}\int_0^t
{ds\over \vert y+X_0(s)\vert^p}\bigg)\Bigg]\nonumber\\
&\le&\omega_d\E\Bigg[\Big(t^{1/p}+\max_{s\le t}\vert X_0(s)\vert\Big)^d
\exp\bigg\{\gamma^{-1}\theta\sup_{y\in\R^d}\int_0^t
{ds\over \vert y+X_0(s)\vert^p}\bigg\}\Bigg]\nonumber\\
&\le&\omega_d\bigg\{\E\Big(t^{1/p}+\max_{s\le t}\vert X_0(s)
\vert\Big)^{d\over 1-\gamma}\bigg\}^{1-\gamma}
\Bigg(\E\exp\bigg\{\gamma^{-2}\theta\sup_{y\in\R^d}\int_0^t {ds\over
\vert y+X_0(s)\vert^p}\bigg\}\Bigg)^\gamma.\nonumber
\end{eqnarray}

By Theorem 1.3 in \cite{BCR},
\begin{equation}\label{eq0422d}
\E\exp\bigg\{\gamma^{-2}\theta\sup_{y\in\R^d}\int_0^t {ds\over \vert
y+X_0(s)\vert^p}\bigg\}<\infty.
\end{equation}
It follows from (\ref{cata-1})-(\ref{eq0422d}) that
$$
\int_{\R^d} \E\Psi\bigg(\theta\int_0^t{ds\over \vert
x+X_0(s)\vert^p}\bigg)dx<\infty.
$$
From (\ref{moment-7}) and (\ref{moment-12}), therefore, we obtain (\ref{th-8}).

Further, according to Theorem 1.3, \cite{BCR},
$$
\begin{aligned}
&\lim_{t\to\infty}{1\over t}\log\E\exp\bigg\{\gamma^{-2}
\theta\sup_{y\in\R^d}\int_0^t
{ds\over \vert y+X_0(s)\vert^p}\bigg\}\\
&=\sup_{g\in {\cal F}_d}\bigg\{\gamma^{-2}\theta\sigma^{-p}
\int_{\R^d}{g^2(x)\over\vert x\vert^p}dx-{1\over 2}
\int_{\R^d}\vert\nabla g(x)\vert^2dx\bigg\}\\
&={2-p\over 2}p^{p\over 2-p}
\Big(\gamma^{-2}\theta\sigma^{-p}\gamma(d,p)\Big)^{2\over 2-p},
\end{aligned}
$$
where ${\cal F}_d=\{g\in W^{1,2}(\R^d);\hskip.1in\| g\|_2=1\}$,
and the last step follows from Lemma 7.3, \cite{Chen-1}.

By (\ref{cata-1}), (\ref{moment-7}) and (\ref{moment-12}) again,
$$
\limsup_{t\to\infty}{1\over t}\log\log
\E_0\otimes\E\exp\bigg\{\theta\int_0^t\ol{V}(s, B_s)ds\bigg\} \le
{2-p\over 2}p^{p\over 2-p}
\Big(\gamma^{-2}\theta\sigma^{-p}\gamma(d,p)\Big)^{2\over 2-p}.
$$
Letting $r\to 1^-$ on the right hand side leads to the desired
upper bound for (\ref{th-9}).

We now establish the lower bound for  (\ref{th-9}). Given
$\epsilon>0$, we estimate as follows
\begin{align}\label{cata-2}
&\E_0\exp\bigg\{\int_{\R^d}\E\Psi\bigg(\theta
\int_0^t{ds\over \vert x+X_0(s)-B_s\vert^p}\bigg)dx\bigg\}\\
&\ge\P_0\Big\{\max_{s\le t}\vert B_s\vert\le {\epsilon\over 2}\Big\}
\exp\bigg\{\int_{\{\vert x\vert\le\frac{\epsilon}{2}\}}\E\Psi\bigg(\theta
\int_0^t\frac{ds}{\(\vert x+X_0(s)\vert+\frac{\epsilon}{2}\)^p}\bigg)dx\bigg\}\nonumber\\
&\ge \P_0\Big\{\max_{s\le t}\vert B_s\vert\le {\epsilon\over
2}\Big\} \exp\bigg\{\omega_d\Big({\epsilon\over
2}\Big)^d\E\Psi\bigg(\theta \int_0^t{ds\over \big(\epsilon+\vert
X_0(s)\vert\big)^p}\bigg)\bigg\}.\nonumber
\end{align}

By the fact that $\Psi(b)\sim e^b$ as $b\to\infty$,
$$
\begin{aligned}
&\liminf_{t\to\infty}{1\over t}\log\E\Psi\bigg(\theta
\int_0^t{ds\over \big(\epsilon+\vert
X_0(s)\vert\big)^p}\bigg)\\
&=\liminf_{t\to\infty}{1\over t}\log\E\exp\bigg\{\theta
\int_0^t{ds\over \big(\epsilon+\vert
X_0(s)\vert\big)^p}\bigg\}\\
&=\sup_{g\in{\cal F}_d}\bigg\{\theta\int_{\R^d} {g^2(x)\over
(\epsilon +\sigma\vert x\vert)^p}dx-{1\over 2}\int_{\R^d}\vert
\nabla g(x)\vert^2dx\bigg\},
\end{aligned}
$$
where the last step follows from a standard treatment of LDP
by Feynman-Kac formula (see, e.g., Theorem 1.6, Chapter 4, \cite{Chen}).

Since $\epsilon>0$ can be arbitrarily small, we obtain
\begin{align}\label{cata-3}
&\liminf_{t\to\infty}{1\over t}\log\log
\E_0\exp\bigg\{\int_{\R^3}\E\Psi\bigg(\theta
\int_0^t{ds\over \vert x+X_0(s)-B_s\vert^p}\bigg)dx\bigg\}\\
&\ge\sup_{g\in{\cal F}_d}\bigg\{\theta\sigma^{-p} \int_{\R^d}
{g^2(x)\over \vert x\vert^p}dx-{1\over 2}\int_{\R^d}\vert \nabla
g(x)\vert^2dx\bigg\}.\nonumber
\end{align}
Thus, the desired lower bound for (\ref{th-9}) follows from
(\ref{moment-7}) and the fact (Lemma 7.3, \cite{Chen-1}) that the
supremum above is equal to the constant appearing on the right hand
side of (\ref{th-9}). \qed

\subsection{Critical case $p=2$ and $d=3$}

For this critical case, a much more delicate treatment than that of
last subsection is needed.

\begin{proposition}
If $\theta>\sigma^2/8$, then
\begin{align}\label{cata-4}
\E_0\exp\bigg\{\int_{\R^3}\E\Psi\bigg(\theta \int_0^t{ds\over \vert
x+X_0(s)-B_s\vert^p}\bigg)dx\bigg\} =\infty,\hskip.2in (t>0).
\end{align}
\end{proposition}
Proof:
Let $t>0$ be fixed.
Given  $\epsilon>0$, by modifying (\ref{cata-2}) slightly,
\begin{align}\label{cata-5}
&\E_0\exp\bigg\{\int_{\R^3}\E\Psi\bigg(\theta
\int_0^t{ds\over \vert x+X_0(s)-B_s\vert^2}\bigg)dx\bigg\}\\
&\ge \P_0\Big\{\max_{s\le t}\vert B_s\vert\le \epsilon\Big\}
\exp\bigg\{\int_{\{\vert x\vert\le 1\}}\E\Psi\bigg(\theta
\int_0^t{ds\over \big(\epsilon +\vert x+X_0(s)\vert\big)^2}
\bigg)dx\bigg\}.\nonumber
\end{align}

Let $0<\delta<t$ be fixed but arbitrary. By (3.16) in Lemma 3.5, \cite{CR}
$$
\begin{aligned}
&\int_{\{\vert x\vert\le 1\}}\E\exp\bigg\{\theta
\int_0^t{ds\over \big(\epsilon +\vert x+X_0(s)\vert\big)^2}
\bigg\}dx\\
&\ge (2\pi\delta)^{3/2}\exp\bigg\{-\delta\epsilon^{-2}+ t\sup_{g\in
{\cal F}_3(B(0, 1))}\bigg(\theta
\int_{B(0,1)}{g^2(x)\over(\epsilon+\sigma\vert x\vert)^2}dx -{1\over
2}\int_{B(0,1)}\vert\nabla g(x)\vert^2dx\bigg)\bigg\}.
\end{aligned}
$$

By Lemma \ref{H-2}, there is a small $a>0$ such that
$$
\sup_{g\in {\cal F}_3(B(0,a\epsilon^{-1})) }
\bigg\{\theta\int_{\{\vert x\vert\le a\epsilon^{-1}\}}
{g^2(x)\over (a+\sigma\vert x\vert)^2}dx-
{1\over 2}\int_{\{\vert x\vert\le a\epsilon^{-1}\}}\vert\nabla g(x)\vert^2dx\bigg\}
\ge t^{-1}
$$
for sufficiently small $\epsilon>0$.

By the substitution $g(x)=\Big({\epsilon\over a}\Big)^{3/2}
f\Big({\epsilon\over a}x\Big)$, therefore,
$$
\sup_{g\in {\cal F}_3(B(0, 1))}\bigg(\theta
\int_{B(0,1)}{g^2(x)\over(\epsilon+\sigma\vert x\vert)^2}dx -{1\over
2}\int_{B(0,1)}\vert\nabla g(x)\vert^2dx\bigg) \ge
t^{-1}\Big({a\over \epsilon}\Big)^2.
$$
Take $\delta<a^2$. We conclude that there is $\gamma >0$ such that
$$
\int_{\{\vert x\vert\le 1\}}\E\exp\bigg\{\theta
\int_0^t{ds\over \big(\epsilon +\vert x+X_0(s)\vert\big)^2}
\bigg\}dx\ge \exp\Big\{\Big({\gamma\over\epsilon}\Big)^2\Big\}
$$
for all small $\epsilon>0$. By the fact that
$$
\begin{aligned}
&\int_{\{\vert x\vert\le 1\}}\E\Psi\bigg(\theta
\int_0^t{ds\over \big(\epsilon +\vert x+X_0(s)\vert\big)^2}
\bigg\}dx\bigg)\\
&=\int_{\{\vert x\vert\le 1\}}\E\exp\bigg\{\theta \int_0^t{ds\over
\big(\epsilon +\vert x+X_0(s)\vert\big)^2} \bigg\}dx -{4\over 3}\pi
-t\int_{\{\vert x\vert\le 1\}}{dx\over \big(\epsilon +\vert
x\vert\big)^2},
\end{aligned}
$$
we have
$$
\int_{\{\vert x\vert\le 1\}}\E\Psi\bigg(\theta \int_0^t{ds\over
\big(\epsilon +\vert x+X_0(s)\vert\big)^2} \bigg\}dx\bigg)\ge
\exp\Big\{\Big({\gamma\over\epsilon}\Big)^2\Big\},
$$
for a possibly different $\gamma>0$.

By the classic estimate
$$
\P_0\Big\{\max_{s\le t}\vert B_s\vert\le \epsilon\Big\} \ge 1-c_1
\exp\Big\{-c_2\epsilon^2\Big\}
$$
for some constants $c_1, \ c_2>0$, and by (\ref{cata-5}),
$$
\begin{aligned}
&\E_0\exp\bigg\{\int_{\R^3}\E\Psi\bigg(\theta
\int_0^t{ds\over \vert x+X_0(s)-B_s\vert^2}\bigg)dx\bigg\}\\
&\ge \left(1-c_1 \exp\Big\{-c_2\epsilon^2\Big\}\right)
\exp\bigg\{\exp\Big\{\Big({\gamma\over\epsilon}\Big)^2\Big\}\bigg\}
\longrightarrow \infty \hskip.2in (\epsilon\to 0^+).
\end{aligned}
$$
This implies (\ref{cata-4}).
\qed

\begin{proposition}
If $\theta<\sigma^2/8$, then
\begin{align}\label{cata-7}
\int_{\R^3}\E\Psi\bigg(\theta\int_0^1{ds\over\vert
x+X_0(s)\vert^2}\bigg) dx<\infty.
\end{align}
\end{proposition}
Proof:
Let $0<\gamma <1$ be sufficiently close to 1 so that
$\gamma^{-1}\theta<\sigma^2/8$. By convexity
$$
\begin{aligned}
&\int_{\R^3}\E\Psi\bigg(\theta\int_0^1{ds\over\vert x+X_0(s)\vert^2}\bigg)dx\\
&\le (1-\gamma)\int_{\R^3}\E\Psi\bigg({\theta\over 1-\gamma}
\int_0^1{1\{\vert x+X_0(s)\vert>1\}\over\vert x+X_0(s)\vert^2}ds\bigg)dx\\
&+\gamma\int_{\R^3}\E\Psi\bigg({\theta\over\gamma} \int_0^1{1\{\vert
x+X_0(s)\vert\le 1\}\over\vert x+X_0(s)\vert^2}ds\bigg)dx.
\end{aligned}
$$
By Jensen's inequality,
$$
\begin{aligned}
&\int_{\R^3}\Psi\bigg({\theta\over 1-\gamma}
\int_0^1{1\{\vert x+X_0(s)\vert>1\}\over\vert x+X_0(s)\vert^2}ds\bigg)dx\\
&\le\int_{\R^3}\int_0^1\Psi\bigg({\theta\over 1-\gamma}
{1\{\vert x+X_0(s)\vert>1\}\over\vert x+X_0(s)\vert^2}\bigg)dxds\\
&=\int_{\{\vert x\vert\ge 1\}} \Psi\Big({\theta\over
1-\gamma}{1\over\vert x\vert^2}\Big)dx<\infty.
\end{aligned}
$$

To prove (\ref{cata-7}), therefore, all we need
is to show is that for any $\theta<\sigma^2/8$,
\begin{align}\label{cata-8}
\int_{\R^3}\E\Psi\bigg(\theta\int_0^1{1\{\vert x+X_0(s)\vert\le 1\}
\over\vert x+X_0(s)\vert^2}ds\bigg) dx<\infty.
\end{align}

Write $\tau_n=\inf\{s\ge 0;\hskip.1in\vert X_0(s)\vert\ge 2^n\}$
($n=0, 1,\cdots$). We have
$$
\begin{aligned}
&\int_{\R^3}
\E\Psi\bigg(\theta\int_0^1{1\{\vert x+X_0(s)\vert\le 1\}\over
\vert x+X_0(s)\vert^2}ds\bigg)dx\\
&=\int_{\R^3}
\E\Bigg[\Psi\bigg(\theta\int_0^1{1\{\vert x+X_0(s)\vert\le 1\}\over
\vert x+X_0(s)\vert^2}ds\bigg);\hskip.05in \tau_0\ge 1\Bigg]dx\\
&+\sum_{n=0}^\infty\int_{\R^3}
\E\Bigg[\Psi\bigg(\theta\int_0^1{1\{\vert x+X_0(s)\vert\le 1\}\over
\vert x+X_0(s)\vert^2}ds\bigg);\hskip.05in \tau_n< 1\le
\tau_{n+1}\Bigg]dx.
\end{aligned}
$$
By the fact that $\Psi(0)=0$,
$$
\begin{aligned}
&\int_{\R^3}
\E\Bigg[\Psi\bigg(\theta\int_0^1{1\{\vert x+X_0(s)\vert\le 1\}\over
\vert x+X_0(s)\vert^2}ds\bigg);\hskip.05in \tau_0\ge 1\Bigg]dx\\
&=\int_{\{\vert x\vert\le 2\}}
\E\Bigg[\Psi\bigg(\theta\int_0^1{1\{\vert x+X_0(s)\vert\le 1\}\over
\vert x+X_0(s)\vert^2}ds\bigg);\hskip.05in \tau_0\ge 1\Bigg]dx
\end{aligned}
$$
and
$$
\begin{aligned}
&\int_{\R^3}
\E\Bigg[\Psi\bigg(\theta\int_0^1{1\{\vert x+X_0(s)\vert\le 1\}\over
\vert x+X_0(s)\vert^2}ds\bigg);\hskip.05in \tau_n< 1\le \tau_{n+1}\Bigg]dx\\
&=\int_{\{\vert x\vert\le 2^{n+1}+1\}}
\E\Bigg[\Psi\bigg(\theta\int_0^1{1\{\vert x+ X_0(s)\vert\le 1\}\over
\vert x+X_0(s)\vert^2}ds\bigg);\hskip.05in \tau_n< 1\le \tau_{n+1}\Bigg]dx\\
&\le\int_{\{\vert x\vert\le 2^{n+2}\}}
\E\Bigg[\Psi\bigg(\theta\int_0^1{1\{\vert x+X_0(s)\vert\le 1\}\over
\vert x+X_0(s)\vert^2}ds\bigg);\hskip.05in \tau_n< 1\le
\tau_{n+1}\Bigg]dx.
\end{aligned}
$$
By the inequality $\Psi(b)\le e^b$ for $b\ge 0$, therefore,
\begin{align}\label{cata-9}
&\int_{\R^3}
\E\Psi\bigg(\theta\int_0^1{1\{\vert x+X_0(s)\vert\le 1\}\over
\vert x+X_0(s)\vert^2}ds\bigg)dx\\
&\le\int_{\{\vert x\vert\le 2\}}
\E\Bigg[\exp\bigg\{\theta\int_0^1{1\{\vert x+X_0(s)\vert\le 1\}\over
\vert x+X_0(s)\vert^2}ds\bigg\}; \tau_0\ge 1\Bigg]dx\nonumber\\
&+\sum_{n=0}^\infty\int_{\{\vert x\vert\le 2^{n+2}\}}\
\E\Bigg[\exp\bigg\{\theta\int_0^1{1\{\vert x+X_0(s)\vert\le 1\}\over
\vert x+X_0(s)\vert^2}ds\bigg\}; \tau_{n}<1\le
\tau_{n+1}\Bigg]dx.\nonumber
\end{align}

Take $\alpha,\beta>1$ such that $\alpha^{-1}+\beta^{-1}=1$ and that
$\alpha \theta<\sigma^2/8$. By H\"older inequality
$$
\begin{aligned}
&\int_{\{\vert x\vert\le 2^{n+2}\}}
\E\Bigg[\exp\bigg\{\theta\int_0^1{1\{\vert x+X_0(s)\vert\le 1\}\over
\vert x+X_0(s)\vert^2}ds\bigg\}; \tau_{n}<1\le \tau_{n+1}\Bigg]dx\\
&\le\Big\{{4\over 3}\pi 2^{3(n+2)}\P\{\tau_{n}<1\}\Big\}^{1/\beta}\\
&\times\Bigg\{\int_{\{\vert x\vert\le 2^{n+2}\}}
\E\Bigg[\exp\bigg\{\alpha\theta\int_0^1{1\{\vert x+X_0(s)\vert\le
1\}\over \vert x+X_0(s)\vert^2}ds\bigg\}; \tau_{n+1}\ge 1\Bigg]dx
\Bigg\}^{1/\alpha}.
\end{aligned}
$$

For each $x\in\R^3$ with $\vert x\vert\le 2^{n+2}$, write
$$
T_x=\inf\Big\{s\ge 0; \hskip.1in \vert x+X_0(s)\vert\ge
2^{n+3}\Big\}.
$$
Then we have
$$
\begin{aligned}
&\int_{\{\vert x\vert\le 2^{n+2}\}}
\E\Bigg[\exp\bigg\{\alpha\theta\int_0^1{1\{\vert x+X_0(s)\vert\le 1\}\over
\vert x+X_0(s)\vert^2}ds\bigg\}; \tau_{n+1}\ge 1\Bigg]dx\\
&\le \int_{\{\vert x\vert\le 2^{n+3}\}}
\E\Bigg[\exp\bigg\{\alpha\theta\int_0^1{ds\over
\vert x+X_0(s)\vert^2}\bigg\}; T_x\ge 1\Bigg]dx\\
&\le {4\over 3}\pi 2^{3(n+3)}\exp\bigg\{\sup_{g\in {\cal F}_d(B(0,
2^{n+3}))} \bigg(\alpha\theta\sigma^{-2}\int_{B(0, 2^{n+3})}{g^2(x)
\over\vert x\vert^2}dx-{1\over 2} \int_{B(0, 2^{n+3})}\vert \nabla
g(x)\vert^2dx\bigg)\bigg\},
\end{aligned}
$$
where the last step follows from Lemma 4.1, \cite{Chen-1}.

Notice that $\alpha\theta\sigma^{-2}<1/8$. By (\ref{H-3}) the $g$-variation
on the right hand side is equal to zero. Hence,
$$
\begin{aligned}
&\int_{\{\vert x\vert\le 2^{n+2}\}}\
\E\Bigg[\exp\bigg\{\theta\int_0^1{1\{\vert x+X_0(s)\vert\le 1\}\over
\vert x+X_0(s)\vert^2}ds\bigg\}; \tau_{n}<1\le \tau_{n+1}\Bigg]dx\\
&\le {4\over 3}\pi 2^{3(n+3)}\Big(\P\{\tau_{n}<1\}\Big)^{1/\beta}
\le C_12^{3n}\exp\Big\{-C_2 2^{2n}\Big\}.
\end{aligned}
$$
A similar argument shows that the first term on the right hand side
of (\ref{cata-9}) is finite. Hence, we have proved (\ref{cata-8}),
and therefore (\ref{cata-7}).

For any $t>0$, by (\ref{moment-12})
\begin{align}\label{cata-10}
&\E_0\exp\bigg\{\int_{\R^3}\E\Psi\bigg(\theta\int_0^t
{ds\over\vert x+X_0(s)-B_s\vert^2}\bigg)dx\bigg\}\\
&\le \exp\bigg\{\int_{\R^3}\E\Psi\bigg(\theta\int_0^t
{ds\over\vert x+X_0(s)\vert^2}\bigg)dx\bigg\}\nonumber\\
&=\exp\bigg\{t^{3/2}\int_{\R^3}\E\Psi\bigg(\theta\int_0^1
{ds\over\vert x+X_0(s)\vert^2}\bigg)dx\bigg\},\nonumber
\end{align}
where the last step follows from Brownian scaling and variable substitution.
By (\ref{cata-7}) the right hand side is finite. Combining this bound with
(\ref{cata-4}), by (\ref{moment-7}) (with $K(x)=\theta\vert x\vert^{-p}$)
 we have proved (\ref{th-10}).
\qed

Finally, we are ready to prove the limit (\ref{th-11}).
Using (\ref{moment-12}) and (\ref{moment-7}), the
upper bound for (\ref{th-11}) follows from (\ref{cata-10}) directly. As for the lower bound, notice that
$$
\begin{aligned}
&\E_0\exp\bigg\{\int_{\R^3}\E\Psi\bigg(\theta\int_0^t
{ds\over\vert x+X_0(s)-B_s\vert^2}\bigg)dx\bigg\}\\
&=\E_0\exp\bigg\{t^{3/2}\int_{\R^3}\E\Psi\bigg(\theta\int_0^1
{ds\over\vert x+X_0(s)-B_s\vert^2}\bigg)dx\bigg\}
\end{aligned}
$$
Thus, the lower bound follows
from the same argument as in Subsection
\ref{R2}.
\qed

\subsection{Super-critical case $p>2$}

Let $t>0$ be fixed. By Jensen's inequality,
$$
\begin{aligned}
&\E_0\exp\bigg\{\int_{\R^d}
\E\Psi\bigg(\theta\int_0^t{ds\over\vert x+X(s)-B_s\vert^p}\bigg)
dx\bigg\}\\
&\ge\exp\bigg\{\int_{\R^d}
\E_0\otimes\E\Psi\bigg(\theta\int_0^t{ds\over\vert x+X_0(s)-B_s\vert^p}\bigg)
dx\bigg\}\\
&=\exp\bigg\{\int_{\R^d}\E\Psi\bigg(\theta\int_0^t {ds\over\vert
x+\sqrt{1+\sigma^{-2}}X_0(s)\vert^p}\bigg)dx\bigg\}.
\end{aligned}
$$
Hence, (\ref{th-12}) follows from (\ref{moment-7}) and the estimate
$$
\begin{aligned}
&\int_{\R^d}\E\Psi\bigg(\theta\int_0^t
{ds\over\vert x+\sqrt{1+\sigma^{-2}}X_0(s)\vert^p}\bigg)dx\\
&\ge \P\Big\{\max_{s\le t}\vert X_0(s)|
\le {\epsilon\over\sqrt{1+\sigma^{-2}}}\Big\}\int_{\{\vert x\vert\ge 2\epsilon\}}
\Psi\Big({t\theta\over 2^p\vert x\vert^p}\Big\}dx\\
&\ge
\exp\big\{-c_1\epsilon^{-2}\big\}\exp\big\{-c_2\epsilon^{-p}\big\}
\longrightarrow \infty,\hskip.2in (\epsilon\to 0^+).
\end{aligned}
$$
\qed

\section{Appendix}\label{a}

In this section, we prove some technical results used in the main body of the paper.

\subsection{Hardy inequality}\label{H}

Hardy's inequality plays an important role in this paper. Searching in literature, we have found large
amount of versions of Hardy's inequality (i.e., \cite{HPL} and \cite{OK}) but the form needed in this paper. For reader's convenience,
we state  Hardy's inequality for $d=3$ in the following lemma and
provide a short proof.

\begin{lemma}\label{H-0'} For any $f\in W^{1,2}(\R^3)$,
\begin{align}\label{H-0}
\int_{\R^3}{f^2(x)\over\vert x\vert^2}dx\le 4\int_{\R^3}\vert\nabla f(x)\vert^2
dx.
\end{align}
Further, the number
4 is the best constant in the sense that for any $\epsilon>0$ one
can find a function $f_\epsilon\in W^{1,2}(\R^3)$ with compact support such that
\begin{align}\label{H-1}
\int_{\R^3}{f_\epsilon^2(x)\over\vert x\vert^2}dx> (4-\epsilon)
\int_{\R^3}\vert\nabla f_\epsilon(x)\vert^2
dx.
\end{align}
\end{lemma}

\proof Write $x=(x_1,x_2,x_3)$. Using integration by parts
$$
\int_{\R^3}{f^2(x)\over\vert x\vert^2}dx
=\int_{\R^3}x_j\Big[{2x_i\over\vert x\vert^4}f^2(x)-{2\over\vert x\vert^2}
f(x){\partial f\over\partial x_j}\Big]dx\hskip.2in j=1,2,3.
$$
Summing over $j$ on the both sides
$$
3\int_{\R^3}{f^2(x)\over\vert x\vert^2}dx=
2\int_{\R^3}\Big[{f^2(x)\over\vert x\vert^2}-
{\nabla f\cdot x\over\vert x\vert^2}f(x)\Big]dx.
$$
Thus,
$$
\int_{\R^3}{f^2(x)\over\vert x\vert^2}dx
=-2\int_{\R^3}{\nabla f\cdot x\over\vert x\vert}{f(x)\over\vert x\vert}dx
\le 2\bigg(\int_{\R^3}{\vert \nabla f\cdot x\vert^2\over\vert x\vert^2}
dx\bigg)^{1/2}\bigg(\int_{\R^3}{f^2(x)\over\vert x\vert^2}dx\bigg)^{1/2}.
$$
Therefore,
$$
\int_{\R^3}{f^2(x)\over\vert x\vert^2}dx\le 4
\int_{\R^3}{\vert \nabla f\cdot x\vert^2\over\vert x\vert^2}dx
\le 4\int_{\R^3}\vert \nabla f(x)\vert^2dx.
$$

To establish (\ref{H-1}), for each large $M>0$, we define
$g_M\in W^{1,2}(\R^3)$
as following:
$$
\begin{aligned}
g_M(x)=\left\{\begin{array}{ll}M^{1/2}\hskip.3in 0\le\vert x\vert\le M^{-1}\\\\
\vert x\vert^{-1/2}\hskip.3in M^{-1}<\vert x\vert\le M\\\\
\displaystyle{2M-\vert x\vert\over M^{3/2}}\hskip.3in M<\vert x\vert\le 2M\\\\
0\hskip.9in\vert x\vert> 2M.\end{array}\right.
\end{aligned}
$$
It is straightforward to exam that $g_M$ is locally supported and
$$
\int_{\R^3}{g_M^2(x)\over\vert x\vert^2}dx=
\bigg\{4-28\Big({7\over 3}+{1\over 2}\log M\Big)^{-1}\bigg\}
\int_{\R^3}\vert\nabla g_M(x)\vert^2dx.
$$
For each $\epsilon>0$, take $M>0$ sufficiently large so
$$
28\Big({7\over 3}+{1\over 2}\log M\Big)^{-1}<\epsilon
$$
and let $f_\epsilon(x)=g_M(x)$.\qed

What has been used in this paper is the following version
of Hardy's inequality.

\begin{lemma}\label{H-2} For any $\theta>0$,
\begin{align}\label{H-3}
\sup_{g\in {\cal F}_3}\bigg\{\theta\int_{\R^3}{g^2(x)\over\vert x\vert^2}dx
-{1\over 2}\int_{\R^3}\vert\nabla g(x)\vert^2dx\bigg\}
= \begin{cases}
 0   & \text{if } \theta\le 1/8, \medskip \\
 \infty   & \text{if } \theta > 1/8.
\end{cases}
\end{align}
\end{lemma}

\proof By Hardy's inequality, the left hand side of (\ref{H-3})
is non-positive when $\theta<1/8$. On the other hand, it is no less than
$$
-{1\over 2}\inf_{g\in {\cal F}_3}\int_{\R^3}\vert\nabla g(x)\vert^2dx
$$
which is equal to zero. Thus, for $\theta\le 1/8$,
$$
\sup_{g\in {\cal F}_3}\bigg\{\theta\int_{\R^3}{g^2(x)\over\vert x\vert^2}dx
-{1\over 2}\int_{\R^3}\vert\nabla g(x)\vert^2dx\bigg\}=0.
$$

Assume $\theta> 1/8$. By the optimality of Hardy's inequality described
in (\ref{H-1}),
$$
H(\theta)\equiv\sup_{g\in {\cal F}_3}
\bigg\{\theta\int_{\R^3}{g^2(x)\over\vert x\vert^2}dx
-{1\over 2}\int_{\R^3}\vert\nabla g(x)\vert^2dx\bigg\}>0.
$$
Given $a>0$, the substitution $g(x)=a^{3/2}f(ax)$ leads to
$H(\theta)=a^2H(\theta)$. So $M(\theta)=\infty$. \qed

\subsection{An auxiliary limit result}

Let
\begin{equation}\label{trap-7}
Q(b)=\int_{\{\vert x\vert\ge b\}}
{1\over\vert x\vert^p}\E{1\{\vert x+U\vert\ge b\}\over \vert x+U\vert^p}dx
\hskip.2in (b\ge 0)
\end{equation}
with $U\sim N(0, I_d)$. In this subsection, we give the limiting behaviors of
$Q(b)$ as $b\to 0+$ and as $b\to\infty$, respectively.

\begin{lemma}\label{lem0422a}
\begin{align}\label{trap-8}
\lim_{b\to 0^+}Q(b)&=\E\int_{\R^d}{1\over\vert x\vert^p}{1\over
\vert x+U\vert^p}dx=C(d,p)\E\vert U\vert^{-(2p-d)}\\
&=C(d,p)2^{d-p-1}d\omega_d(2\pi)^{-d/2}\Gamma (d-p),\nonumber
\end{align}
where
\begin{align}\label{trap-2}
C(d,p)=\pi^{d/2}{\displaystyle \Gamma^2\Big({d-p\over 2}\Big)
\Gamma\Big({2p-d\over 2}\Big)\over\displaystyle \Gamma^2\Big({p\over
2}\Big)\Gamma(d-p)}.
\end{align}
\end{lemma}
Proof: The first equality follows from the monotone convergence theorem. The second follows from  the
identity (see p.118, \cite{Donoghue} or p.118, (8), \cite{Stein})
\begin{align}\label{trap-1}
\int_{\R^d}{1\over \vert x-y\vert^p}{1\over \vert x-z\vert^p}dx
=C(d,p){1\over \vert y-z\vert^{2p-d}}\hskip.2in y,z\in\R^d.
\end{align}
\qed

\begin{lemma}\label{a-1}
Under $d/2<p<d$, we have
\begin{align}\label{a-2}
\lim_{b\to\infty}b^{2p-d}Q(b)={d\omega_d\over 2p-d}.
\end{align}
\end{lemma}

\proof By variable substitution, we have,
$$
\begin{aligned}
&\int_{\{\vert x\vert\ge b\}}
{1\over\vert x\vert^p}{1\{\vert x+ U\vert\ge b\}\over\vert
x+U\vert^p}dx\\
&=\vert U\vert^{-(2p-d)}\int_{\{\vert x\vert\ge b\vert U\vert^{-1}\}}
{1\over\vert x\vert^p}{1\{\vert x+\vert U\vert^{-1}U\vert\ge
b\vert U\vert^{-1}\}\over\vert x+\vert U\vert^{-1}U\vert^p}dx\\
&=\vert U\vert^{-(2p-d)}H\big(b\vert U\vert^{-1}\big),
\end{aligned}
$$
where
$$
H(b)=\int_{\{\vert x\vert\ge b\}}{1\over\vert x\vert^p} {1\{\vert
x+x_0\vert\ge b\}\over\vert x+x_0\vert^p}dx,
$$
$x_0$ is a fixed point with $\vert x_0\vert =1$, and the last step
follows from the fact that $H(b)$ does not depend on the location of
$x_0$ on the unit sphere. Thus,
$$
\begin{aligned}
&Q(b)=\E\vert U\vert^{-(2p-d)}H\big(b\vert U\vert^{-1}\big)\\
&=\int_{\R^d}{dx\over\vert x\vert^p\vert x+x_0\vert^p}
\bigg[\int_{\R^d}{1\{\vert y\vert\ge b\vert x\vert^{-1}\}
1\{\vert y\vert\ge b\vert x+x_0\vert^{-1}\}\over \vert y\vert^{2p-d}} p_1(y)
dy\bigg]dx\\
&\sim\int_{\{\vert x\vert\ge C\}}{dx\over\vert x\vert^p\vert
x+x_0\vert^p} \bigg[\int_{\R^d}{1\{\vert y\vert\ge b\vert
x\vert^{-1}\} 1\{\vert y\vert\ge b\vert x+x_0\vert^{-1}\}\over \vert
y\vert^{2p-d}} p_1(y)dy\bigg]dx,
\end{aligned}
$$
where $p_1(x)$ is the density of the $d$-dimensional standard normal
distribution, $C>0$ is a large but fixed constant. By the fact that
$C\gg 1=\vert x_0\vert$, for $b\to\infty$, we have
$$
\begin{aligned}
&Q(b)\sim q(C)\int_{\{\vert x\vert\ge C\}}
{dx\over\vert x\vert^{2p}}
\bigg[\int_{\{\vert y\vert\ge b\vert x\vert^{-1}\}}
{1\over \vert y\vert^{2p-d}}p_1(y)dy\bigg]\\
&=q(C)\int_{\R^d}{1\over \vert y\vert^{2p-d}}p_1(y)\bigg[
\int_{\{\vert x\vert\ge
\max\{C, b\vert y\vert^{-1}\}\}}{dx\over\vert x\vert^{2p}}\bigg]
dy\\
&=q(C){d\omega_d\over 2p-d}\int_{\R^d}{1\over \vert y\vert^{2p-d}}p_1(y)
\min\bigg\{{1\over C^{2p-d}}, \Big({\vert y\vert\over b}\Big)^{2p-d}\bigg\}dy\\
&=q(C){d\omega_d\over 2p-d}{1\over C^{2p-d}}\int_{\{\vert y\vert\ge C^{-1}b\}}
{1\over \vert y\vert^{2p-d}}p_1(y)dy\\
&+q(C){d\omega_d\over 2p-d}{1\over b^{2p-d}} \int_{\{\vert y\vert\le
C^{-1}b\}}p_1(y)dy,
\end{aligned}
$$
where $q(C)\to 1$ as $C\to\infty$.
The first term on the right hand side is obviously negligible. \qed

\section*{Acknowledgment}

The authors would like to thank the anonymous referee for helpful suggestions.

\bibliographystyle{amsplain}

\end{document}